\documentclass[11pt, twoside, leqno]{article}
\usepackage{amssymb,amsmath,amsfonts,amsthm,mathrsfs}
\usepackage{indentfirst}
\usepackage{graphics}
\usepackage{color}
\usepackage[margin=3cm]{geometry}
\usepackage{enumerate}
\usepackage{amsmath}
\usepackage{bbm}           % used for blackboard 1
\usepackage{bm}
\usepackage{eso-pic,graphicx}
\usepackage{tikz}
\usepackage{cite}
\usepackage{esint}
\usepackage[colorlinks=true, pdfstartview=FitV, linkcolor=cyan, citecolor=magenta, urlcolor=magenta]{hyperref}
\usepackage{pslatex}
\usepackage{rotating}
\usepackage{verbatim}

\makeatletter
\def\Ddots{\mathinner{\mkern1mu\raise\p@
		\vbox{\kern7\p@\hbox{.}}\mkern2mu
		\raise4\p@\hbox{.}\mkern2mu\raise7\p@\hbox{.}\mkern1mu}}
\makeatother

\def\XXint#1#2#3{{\setbox0=\hbox{$#1{#2#3}{\int}$}
		\vcenter{\hbox{$#2#3$}}\kern-.5\wd0}}

\begin{document}
\def\rn{{\mathbb R^n}}  \def\sn{{\mathbb S^{n-1}}}
\def\co{{\mathcal C_\Omega}}
\def\z{{\mathbb Z}}
\def\nm{{\mathbb (\rn)^m}}
\def\mm{{\mathbb (\rn)^{m+1}}}
\def\n{{\mathbb N}}
\def\cc{{\mathbb C}}

\newtheorem{defn}{Definition}
\newtheorem{thm}{Theorem}
\newtheorem{lem}{Lemma}
\newtheorem{cor}{Corollary}
\newtheorem{rem}{Remark}
\renewcommand{\theequation}{\arabic{section}.\arabic{equation}}

\title{\bf\Large The multilinear Littlewood-Paley square operators and their commutators on weighted Morrey spaces
\footnotetext{{\it Key words and phrases}: multilinear square operator; Littlewood-Paley operator; weighted Morrey space; commutator; weak-type $L\log L$ estimate; multilinear Marcinkiewicz integral; BMO.
\newline\indent\hspace{1mm} {\it 2020 Mathematics Subject Classification}: Primary 42B25; Secondary 42B20, 47H60, 47B47.
\newline\indent\hspace{1mm} E-mail: xicenmath@gmail.com.
\newline\indent\hspace{1mm} Address: School of Mathematics and Physics, Southwest University of Science and Technology, Mianyang, 621010, P. R. China.}}
\date{}
\author{Xi Cen}
\maketitle

\begin{center}
\begin{minipage}{13cm}
{\small {\bf Abstract:}\quad
In this paper, we prove the boundedness of the multilinear Littlewood-Paley square operators and their commutators on weighted Morrey spaces, then we give the boundedness and weak-type $L\log L$ estimates for the commutators of multilinear Littlewood-Paley g-functions and multilinear Marcinkiewicz integrals on weighted Morrey spaces in the form of corollaries.}
\end{minipage}
\end{center}

%------------------------section 1------------------------
\section{Introduction}\label{sec1}
\subsection{Background}
\par
It is well known that the Littlewood-Paley $g$-function is a very important tool in harmonic analysis and the Marcinkiewicz integral is essentially a Littlewood-Paley $g$-function. The Littlewood-Paley $g$-function in one dimension was first introduced by Littlewood and Paley in studying the dyadic decomposition of Fourier series, and this theory is extended to higher dimension by Stein. In \cite{Wang1,Wang2}, Wang had proved the boundedness of multilinear Calder\'on--Zygmund operators, multilinear fractional integrals and Marcinkiewicz integrals with rough kernel on the weighted Morrey spaces $L^{p,\kappa}(w).$ The singular integral operators and Littlewood-Paley $g$-functions play an important role in partial differential equations, so many mathematicians have studied them in different spaces, see \cite{Hormozi,Lerner,LuYangZhou,ShiFu,xue1,xue2,xue3,xue4,xue5,xue6,Wang3,Wang2,Wang1}. 

\par
%定义1
\begin{defn}
Suppose that $\varphi  \in L_{loc}^1(\rn)$, the Littlewood-Paley $g$-function is defined by$${g_\varphi }(f) = {({\int_0^{ \infty } {\left| {f * {\varphi _t}} \right|} ^2}\frac{{dt}}{t})^{\frac{1}{2}}}.$$	
\end{defn}
From 2013 to 2018, Xue{\cite{xue1,xue2,xue3,xue4,xue5}} studied a multilinear version of this operator, who generalized it to the cases of non-convolution kernels, Dini’s type kernels and non-smooth kernels, and proved the boundedness in m-fold weighted Lebesgue spaces. Now, we cite some definitions as follows.
%定义
\begin{defn}[\cite{xue2}]
Suppose that $w :[0, + \infty ) \to [0, + \infty )$ is a nondecreasing function with $0 < w (1) < \infty.$ For $a>0,$ we say $w  \in Dini(a),$ if $${\left[ w  \right]_{Dini(a)}} = \int_0^1 {\frac{{{w ^a}(t)}}{t}} dt < \infty. $$
\end{defn}

%定义
\begin{defn}[\cite{xue2}]
	For any $t \in (0,\infty ),$ let ${K}(x,{y_1}, \cdots ,{y_m})$ be a locally integrable function defined away from the diagonal $x = {y_1} =  \cdots  = {y_m}$ in $(\mathbb R^n)^{m+1}$ and denote $(x,\vec y) = (x,{y_1}, \cdots ,{y_m})$, ${K_t}(x,\vec y) = \frac{1}{{{t^{mn}}}}{K}(\frac{x}{t},\frac{{{y_1}}}{t}, \ldots ,\frac{{{y_m}}}{t})$. we will always use this notation throughout this paper. We say $K$ is a kernel of type $w$, if for some constants $A>0,$ the following inequalities hold:
	\begin{equation}\label{NW5}
	{(\int_0^\infty  {{{\left| {{K_t}(x,\vec y)} \right|}^2}\frac{{dt}}{t}} )^{\frac{1}{2}}} \le \frac{A}{{{{(\sum\limits_{j = 1}^m {\left| {x - {y_j}} \right|} )}^{mn}}}}
	\end{equation}
	\begin{equation}\label{NW6}
	{(\int_0^\infty  {{{\left| {{K_t}(z,\vec y) - {K_t}(x,\vec y)} \right|}^2}} \frac{{dt}}{t})^{\frac{1}{2}}} \le \frac{A}{{{{(\sum\limits_{j = 1}^m {\left| {x - {y_j}} \right|} )}^{mn}}}} \cdot w (\frac{{\left| {z - x} \right|}}{{\sum\limits_{j = 1}^m {\left| {x - {y_j}} \right|} }}),
	\end{equation}
	whenever $\left| {z - x} \right| \le \frac{1}{2}\mathop {\max }\limits_{1 \le j \le m} \{ \left| {x - {y_j}} \right|\} ;$ and
	\begin{equation}\label{NW7}
	{(\int_0^\infty  {{{\left| {{K_t}(x,\vec y) - {K_t}(x,{y_1}, \cdots ,{{y_i}'}, \cdots ,{y_m})} \right|}^2}} \frac{{dt}}{t})^{\frac{1}{2}}} \le \frac{A}{{{{(\sum\limits_{j = 1}^m {\left| {x - {y_j}} \right|} )}^{mn}}}} \cdot w (\frac{{{{\left| {{y_i} - {{y'}_i}} \right|} }}}{{\sum\limits_{j = 1}^m {\left| {x - {y_j}} \right|} }}),
	\end{equation}
	for any $i \in \{ 1, \cdots ,m\},$ whenever $\left| {{y_i} - {y_i}^\prime } \right| \le \frac{1}{2}\mathop {\max }\limits_{1 \le j \le m} \{ \left| {x - {y_j}} \right|\}$.
\end{defn}
When $w\left( t \right) = {t^\gamma }$ for some $\gamma  > 0$, we say $K$ satisfies the integral condition of C-Z type $I$, which is introduced as follows.
%定义
{\begin{defn}[\cite{xue1}]
For any $t \in (0,\infty ),$ let ${K}(x,{y_1}, \cdots ,{y_m})$ be a locally integrable function defined away from the diagonal $x = {y_1} =  \cdots  = {y_m}$ in $(\mathbb R^n)^{m+1}.$ We say $K$ satisfies the integral condition of C-Z type I, if for some positive constants $\gamma, A,$ and $B>1,$ the following inequalities hold:
\begin{equation}\label{EQQ001}
{(\int_0^\infty  {{{\left| {{K_t}(x,\vec y)} \right|}^2}\frac{{dt}}{t}} )^{\frac{1}{2}}} \le \frac{A}{{{{(\sum\limits_{j = 1}^m {\left| {x - {y_j}} \right|} )}^{mn}}}}
\end{equation}
\begin{equation}\label{EQQ002}
{(\int_0^\infty  {{{\left| {{K_t}(z,\vec y) - {K_t}(x,\vec y)} \right|}^2}} \frac{{dt}}{t})^{\frac{1}{2}}} \le \frac{{A{{\left| {z - x} \right|}^\gamma }}}{{{{(\sum\limits_{j = 1}^m {\left| {x - {y_j}} \right|} )}^{mn + \gamma }}}},
\end{equation}
whenever $\left| {z - x} \right| \le \frac{1}{B}\mathop {\max }\limits_{1 \le j \le m} \{ \left| {x - {y_j}} \right|\} ;$ and
\begin{equation}\label{EQQ003}
{(\int_0^\infty  {{{\left| {{K_t}(x,\vec y) - {K_t}(x,{y_1}, \cdots ,{{y'}_i}, \cdots ,{y_m})} \right|}^2}} \frac{{dt}}{t})^{\frac{1}{2}}} \le \frac{{A{{\left| {{y_i} - {{y'}_i}} \right|}^\gamma }}}{{{{(\sum\limits_{j = 1}^m {\left| {x - {y_j}} \right|} )}^{mn + \gamma }}}},
\end{equation}
for any $i \in \{ 1, \cdots ,m\},$ whenever $\left| {{y_i} - {y_i}^\prime } \right| \le \frac{{\left| {x - {y_i}} \right|}}{B}.$\\
\end{defn}
\begin{defn}[\cite{xue1}]
	For any $t \in (0,\infty ),$ let ${K}(x,{y_1}, \cdots ,{y_m})$ be a locally integrable function defined away from the diagonal $x = {y_1} =  \cdots  = {y_m}$ in $(\mathbb R^n)^{m+1}.$ We say $K$ satisfies the integral condition of C-Z type II, if for some positive constants $\gamma, A,$ and $B>1,$ the following inequalities hold:
	\begin{equation}\label{EQQ005}
	{({\int_0^\infty  {\int_{\rn}^{} {\left| {{{(\frac{t}{{\left| {x - z} \right| + t}})}^{\frac{{n\lambda }}{2}}}{K_t}(z,\vec y)} \right|} } ^2}\frac{{dzdt}}{{{t^{n + 1}}}})^{\frac{1}{2}}} \le \frac{A}{{{{(\sum\limits_{j = 1}^m {\left| {x - {y_j}} \right|} )}^{mn}}}}
	\end{equation}
	\begin{equation}\label{EQQ006}
	(\int_0^\infty  {{{\int_{\rn}^{} {{{(\frac{t}{{\left| z \right| + t}})}^{n\lambda }}\left| {{K_t}(x - z,\vec y) - {K_t}(x' - z,\vec y)} \right|} }^2}\frac{{dzdt}}{{{t^{n + 1}}}}{)^{\frac{1}{2}}}}  \le \frac{{A{{\left| {x - x'} \right|}^\gamma }}}{{{{(\sum\limits_{j = 1}^m {\left| {x - {y_j}} \right|} )}^{mn + \gamma }}}}
	\end{equation}
	whenever $\left| {x - x'} \right| \le \frac{1}{B}\mathop {\max }\limits_{1 \le j \le m} \{ \left| {x - {y_j}} \right|\};$and
	\begin{equation}\label{EQQ007}
	(\int_0^\infty  {{{\int_{\rn}^{} {{{(\frac{t}{{\left| {x - z} \right| + t}})}^{n\lambda }}\left| {{K_t}(z,\vec y) - {K_t}(z,{y_1}, \cdots ,{y_i}^\prime , \cdots ,{y_m})} \right|} }^2}\frac{{dzdt}}{{{t^{n + 1}}}}{)^{\frac{1}{2}}}}  \le \frac{{A{{\left| {{y_i} - {{y'}_i}} \right|}^\gamma }}}{{{{(\sum\limits_{j = 1}^m {\left| {x - {y_j}} \right|} )}^{mn + \gamma }}}}
	\end{equation}
	for $i \in \{ 1, \cdots ,m\},$ whenever $\left| {{y_i} - {y_i}^\prime } \right| \le \frac{{\left| {x - {y_i}} \right|}}{B}$
\end{defn}
\begin{defn}[\cite{xue1}]
	Let $K$ be a function defined on $\mathbb{R}^n\times \mathbb{R}^{mn}$ with $supp K\subseteq
	\mathcal{B}:=\{(x,y_1,\dots,y_m):\sum_{j=1}^m|x-y_j|^2\leq 1\}$. $K$ is called a multilinear Marcinkiewicz kernel if for some $0<\delta<mn$ and some positive constants $A$, $\gamma_0$, and $B_1$,
	\begin{enumerate}
		\item[\emph{(a)}]$
		|K(x,\vec{y})|\leq\frac{A}
		{(\sum_{j=1}^{m}|x-y_{j}|)^{mn-\delta}};
		$	\item[\emph{(b)}]$
		|K(x,\vec{y})-K(x,y_{1},\dots,y_i',\dots,y_{m})|
		\leq\frac{A|y_i-y_i'|^{\gamma_0}}
		{(\sum_{j=1}^{m}|x-y_{j}|)^{mn-\delta+\gamma_0}},
		$;
		\item[\emph{(c)}]$
		|K(x,\vec{y})-K(x',y_{1},\dots,y_{m})|
		\leq\frac{A|x-x'|^{\gamma_0}}
		{(\sum_{j=1}^{m}|x-y_{j}|)^{mn-\delta+\gamma_0}},
		$,\end{enumerate}
	where (b) holds whenever $(x,y_1,\dots,y_m)\in \mathcal{B}$ and
	$|y_i-y_i'|\leq\frac{1}{B_1}|x-y_i|$ for all $0\leq i\leq m$,
	and (c) holds whenever $(x,y_1,\dots,y_m)\in \mathcal{B}$ and
	$|x-x'|\leq\frac{1}{B_1}\max_{1\leq j \leq m}|x-y_{j}|$.
\end{defn}

\begin{defn}[\cite{xue1}]
	Let $K(x,y_1,\dots,y_m)$ be a locally integrable function defined away from the diagonal
	$x=y_1=\dots=y_m$ in $(\mathbb{R}^n)^{m+1}$.
	$K$ is called a multilinear Littlewood-Paley kernel if for some positive constants $A$, $\gamma_0$, $\delta$, and $B_1$, it holds that 
	\begin{enumerate}
		\item[\emph{(d)}]$
		|K(x,\vec{y})|\leq\frac{A}
		{(1+\sum_{j=1}^{m}|x-y_{j}|)^{mn+\delta}};
		$\item[\emph{(e)}]$
		|K(x,\vec{y})-K(x,y_{1},\dots,y_i',\dots,y_{m})|
		\leq\frac{A|y_i-y_i'|^{\gamma_0}}
		{(1+\sum_{j=1}^{m}|x-y_{j}|)^{mn+\delta+\gamma_0}}
		$;
		\item[\emph{(f)}]$
		|K(x,\vec{y})-K(x',y_{1},\dots,y_{m})|
		\leq\frac{A|x-x'|^{\gamma_0}}
		{(1+\sum_{j=1}^{m}|x-y_{j}|)^{mn+\delta+\gamma_0}},
		$\end{enumerate}
	where (e) holds whenever
	$|y_i-y_i'|\leq\frac{1}{B_1}|x-y_i|$ and for all $1\leq i\leq m$,
	and (f) holds whenever
	$|x-x'|\leq\frac{1}{B_1}\max\limits_{1\leq j \leq m}|x-y_{j}|$.
\end{defn}
The following two lemmas are crucial for understanding.
\begin{lem}[\cite{xue1}]
	If $K$ is either a multilinear Littlewood-Paley kernel or multilinear Marcinkiewicz kernel, then $K$ is a kernel of C-Z  type I.
\end{lem}
\begin{lem}[\cite{xue1}]
	If $K$ is multilinear Littlewood-Paley kernel, then $K$ is a kernel of C-Z type II.
\end{lem}

The multilinear square function with a kernel of C-Z type $I$ or with kernel of type $w(t),$ and $w \in Dini(1)$ is defined by
\begin{equation}\label{NW8}
T(\vec f)(x) = {({\int_0^\infty  {\left| {\int_{\nm} {{K_t}(x,\vec y)\prod\limits_{j = 1}^m {{f_j}({y_j})d{y_1} \cdots d{y_m}} } } \right|} ^2}\frac{{dt}}{t})^{\frac{1}{2}}},
\end{equation}
for any $\vec f = ({f_1}, \cdots ,{f_m}) \in {\mathscr{S}}({\rm{\rn}}) \times  \cdots  \times {\mathscr{S}}({\rm{\rn}})$ and all $x \notin \bigcap\limits_{j = 1}^m {{\rm{supp}}{f_j}},$
and assume that T can be extended to be a bounded operator from ${L^{{q_1}}} \times  \cdots {L^{{q_m}}}$ to ${L^q}$, for some $1 \le {q_1} \cdots , {q_m} \le \infty, \frac{1}{q} = \sum\limits_{k = 1}^m {\frac{1}{{{q_k}}}}$.

$T$ is called a multilinear Marcinkiewicz operator when $K$ is a  multilinear Marcinkiewicz kernel. $T$ is called a multilinear Littlewood-Paley $g$-function when $K$ is a multilinear Littlewood-Paley kernel.

The multilinear square function with a kernel of C-Z type $II$ is defined by
\begin{equation*}\aligned
T_{\lambda}(\vec{f})(x)&=\bigg(\iint_{\mathbb{R}^{n+1}_+}\big(\frac{t}{|x-z|+t}\big)^{n\lambda}
|\int_{\mathbb{R}^{nm}}K_t(z,\vec{y})
\prod_{j=1}^mf_j(y_j)d\vec{y}|^2\frac{dzdt}{t^{n+1}}\bigg)^{\frac12},
\endaligned
\end{equation*}
for any $\vec f = ({f_1}, \cdots ,{f_m}) \in {\mathscr{S}}({\rm{\rn}}) \times  \cdots  \times {\mathscr{S}}({\rm{\rn}})$ and all $x \notin \bigcap\limits_{j = 1}^m {{\rm{supp}}{f_j}},$
and assume that T can be extended to be a bounded operator from ${L^{{q_1}}} \times  \cdots {L^{{q_m}}}$ to ${L^q}$, for some $1 \le {q_1} \cdots , {q_m} \le \infty, \frac{1}{q} = \sum\limits_{k = 1}^m {\frac{1}{{{q_k}}}}$.

$T_{\lambda}$ is called a multilinear Littlewood-Paley $g_\lambda ^*$-function when $K$ is a multilinear Littlewood-Paley kernel.
\par
The classical Morrey spaces $\mathcal L^{p,\lambda}$ were first introduced by Morrey in \cite{Morrey} to study the local behavior of solutions to second order elliptic partial differential equations. In 2009, Komori and Shirai \cite{Komori} considered the weighted version of Morrey spaces $L^{p,\kappa}(\omega)$ and studied the boundedness of some classical operators such as the Hardy-Littlewood maximal operator and the Calder\'on-Zygmund operator on these spaces.
\begin{defn}[\cite{Komori}]
Let $0<p<\infty$, $0<\kappa<1$ and $\omega$ be a weight function on $\mathbb R^n$. Then the weighted Morrey space is defined by
\begin{equation*}
L^{p,\kappa}(\omega)=\big\{f:\big\|f\big\|_{L^{p,\kappa}(\omega)}:= \mathop {\sup }\limits_B \omega {(B)^{ - \frac{\kappa }{p}}}{\left\| f \right\|_{{L^p}(B,\omega dx)}}<\infty\big\}.
\end{equation*}
\end{defn}

\begin{defn}[\cite{Wang1}]
Let $0<p<\infty$, $0<\kappa<1$ and $\omega$ be a weight function on $\mathbb R^n$. Then the weighted weak Morrey space is defined by
\begin{equation*}
WL^{p,\kappa}(\omega)=\big\{f:\big\|f\big\|_{WL^{p,\kappa}(\omega)}:= \mathop {\sup }\limits_B \omega {(B)^{ - \frac{\kappa }{p}}}{\left\| f \right\|_{W{L^p}(B,\omega dx)}}<\infty\big\}.
\end{equation*}
\end{defn}
In order to deal with the end-point case of the commutators, we have to consider the following $L\log L$-type space.
\begin{defn}[\cite{Wang4}]
	Let $p=1$, $0<\kappa<1$ and let $\omega$ be a weight on $\mathbb R^n$. We denote by $(L\log L)^{1,\kappa}(\omega)$ the weighted Morrey space of $L\log L$ type, which is defined by
	\begin{equation*}
	(L\log L)^{1,\kappa}(\omega):=\Big\{f:	\big\|f\big\|_{(L\log L)^{1,\kappa}(\omega)}:=\sup_{B}\omega(B)^{1-\kappa}\big\|f\big\|_{L\log L(\omega),B}<\infty\Big\}.
	\end{equation*}
\end{defn}
Here $\|\cdot\|_{L\log L(\omega),B}$ denotes the weighted Luxemburg norm, whose definition will be given in Section \ref{sec2}. Note that $t\leq t(1+\log^+t)$. For any ball $B$ in $\mathbb R^n$ and $\omega\in A_\infty$, we have a importamt inequalities as follows:
\begin{equation}\label{main esti1}
\big\|f\big\|_{L(\omega),B}=\frac{1}{\omega(B)}\int_{B}|f(x)|\omega(x)\,dx\leq\big\|f\big\|_{L\log L(\omega),B}.
\end{equation}
In fact, for every $\sigma \in E:=\bigg\{\sigma>0:\frac{1}{\omega(B)}\int_B\Phi\bigg(\frac{|f(x)|}{\sigma}\bigg) \omega(x)\,dx\leq1\bigg\}$, we have
$$\frac{1}{{\omega (B)}}\int_B {\frac{{\left| {f(x)} \right|}}{\sigma }\omega (x){\mkern 1mu} dx}  \le \frac{1}{{\omega (B)}}\int_B {\Phi (\frac{{\left| {f(x)} \right|}}{\sigma })\omega (x){\mkern 1mu} dx};$$
then we deduce,
$${\left\| f \right\|_{L(\omega ),B}} \le \mathop {\inf }\limits_{\sigma  \in E} \sigma \frac{1}{{\omega (B)}}\int_B {\Phi (\frac{{\left| {f(x)} \right|}}{\sigma })\omega (x){\mkern 1mu} dx}  \le \mathop {\inf }\limits_{\sigma  \in E} \sigma  = {\left\| f \right\|_{\Phi (\omega ),B}},$$
where $\Phi (t)=t(1+\log^+t)$. Thus we obtain the estimate (\ref{main esti1}).

Many people have studied different types of singular integral operators on weighted Morrey spaces and we present some of their works below.

In 2013, Wang and Yi \cite{Wang1} have proved the boundedness of multilinear Calderón-Zygmund and fractional integral operators on weighted Morrey spaces

In 2014, Iida \cite{Iida} has studied the boundedness of the Hardy-Littlewood maximal operator and multilinear maximal operator in weighted Morrey type spaces. He, Zheng and Tao \cite{He-Zheng-Tao} have obtained the estimates for multilinear commutators of generalized fractional integral operators on weighted Morrey Space. Hu and Wang \cite{Hu-Wang} showed the bounedness of multilinear fractional integral operators on generalized weighted Morrey spaces. Hu, Li and Wang \cite{Hu-Li-Wang} have proved the boundedness of multilinear singular integral operators on generalized weighted Morrey spaces.

In 2016, He and Tao \cite{He-Tao} have established the theory of multilinear singular operators with rough kernels on the weighted Morrey spaces.He and Zhou \cite{He-Zhou} have come up with the boundedness of vector-valued maximal multilinear Calderón–Zygmund operator with nonsmooth kernel on weighted Morrey spaces.

In 2021, Ismayilova \cite{Ismayilova} studied Calderón-Zygmund operators with kernels of Dini’s type and their multilinear commutators on generalized Morrey spaces. Lin and Yan \cite{Lin-Yan} have proved the boundedness of multilinear strongly singular Calderón-Zygmund operators and commutators on Morrey type spaces.

Now, we introduce the main results of this paper.
\subsection{Main Results}
\par
Firstly, we give the boundedness of the multilinear vector-valued operators on weighted Morrey spaces.
%多线性向量值函数
\begin{thm}\label{the1}
Let $m\in\n,$ $X$ is a Banach space, $B(\cc,X)$ is the space of all bounded linear operators from $\cc$ to $X$, suppose that operator-valued function $Q:(\mm \backslash E) \to B(\cc,X),$ $E = \{ (x,\vec y) \in \mm:x = {y_1} =  \cdots  = {y_m}\},$ which satisfies Size Condition:
\begin{equation}
{\left\| {Q(x,\vec y)} \right\|_{B(\cc,X)}} \le \frac{C}{{{{(\sum\limits_{j = 1}^m {\left| {x - {y_j}} \right|} )}^{mn}}}}.
\end{equation}
We define the multilinear $X$-valued operator $T$ by
\begin{equation}
T(\vec f)(x) = {\left\| {\int_{\nm} {(Q(x,\vec y))(\prod\limits_{j = 1}^m {{f_j}({y_j}))d{y_1} \cdots d{y_m}} } } \right\|_X}
\end{equation}
If $p_1,\ldots,p_m\in[1,\infty)$ with $1/p=\sum_{k=1}^m 1/{p_k}$, and $\vec{\omega}=(\omega_1,\ldots,\omega_m)\in {A_{\vec P}} \cap {\left( {{A_\infty }}\right)^m}$. For any $0<\kappa<1$, the following results hold:
\begin{enumerate}[(i)]
\item If $\mathop {\min }\limits_{1 \le i \le m} \{ {p_i}\}  > 1$, such that T is well-defined on ${L^{{p_1}}}({\omega _1}) \times  \cdots  \times {L^{{p_m}}}({\omega _m})$, which is also bounded from ${L^{{p_1}}}({\omega _1}) \times  \cdots  \times {L^{{p_m}}}({\omega _m})$ to ${L^p}({v_{\vec \omega }}),$ then there exists a constant ${C}$, independent of $\vec f$, such that 
\begin{equation*}
{\left\| {T(\vec f)} \right\|_{{L^{p,\kappa}}({v_{\vec \omega }})}} \le C{\prod\limits_{i = 1}^m {\left\| {{f_i}} \right\|} _{{L^{{p_i},\kappa }}({\omega _i})}}.
\end{equation*}
\item If $\mathop {\min }\limits_{1 \le i \le m} \{ {p_i}\}  = 1,$ such that T is well-defined on ${L^{{p_1}}}({\omega _1}) \times  \cdots  \times {L^{{p_m}}}({\omega _m})$, which is also bounded from ${L^{{p_1}}}({\omega _1}) \times  \cdots  \times {L^{{p_m}}}({\omega _m})$ to $W{L^p}({v_{\vec \omega }}),$ then there exists a constant ${C}$, independent of $\vec f$, such that  
\begin{equation*}
{\left\| {T(\vec f)} \right\|_{W{L^{p,\kappa}}({v_{\vec \omega }})}} \le C{\prod\limits_{i = 1}^m {\left\| {{f_i}} \right\|} _{{L^{{p_i},\kappa}}({\omega _i})}}.
\end{equation*}
\end{enumerate}
\end{thm}
Now, according to theorem \ref{the1} and Theorem 3.1 in \cite{xue6}, we have a specific corollary as follows.
\begin{cor}
Let $m\in\n,$ $X$ is a Banach space, $B(\cc,X)$ is the space of all bounded linear operators from $\cc$ to $X$, suppose that operator-valued function $Q:(\mm \backslash E) \to B(\cc,X),$ $E = \{ (x,\vec y) \in \mm:x = {y_1} =  \cdots  = {y_m}\}$, we define the  multilinear $X$-valued Calder\'on--Zygmund operator $T$ by
\begin{equation}
T(\vec f)(x) = {\left\| {\int_{\nm} {(Q(x,\vec y))(\prod\limits_{j = 1}^m {{f_j}({y_j}))d{y_1} \cdots d{y_m}} } } \right\|_X}
\end{equation}
for any $\vec f = ({f_1}, \cdots ,{f_m}) \in {\mathscr{S}}({\rm{\rn}}) \times  \cdots  \times {\mathscr{S}}({\rm{\rn}})$ and all $x \notin \bigcap\limits_{j = 1}^m {{\rm{supp}}{f_j}},$
and assume that T can be extended to be a bounded operator from ${L^{{q_1}}} \times  \cdots {L^{{q_m}}}$ to ${L^q}$, for some $1 \le {q_1} \cdots , {q_m} \le \infty, \frac{1}{q} = \sum\limits_{k = 1}^m {\frac{1}{{{q_k}}}}>0$. The kernel satisfies, for some $\varepsilon,C>0$,
\begin{enumerate}
	\item[\emph{(i)}]${\left\| {Q(x,{y_1}, \ldots ,{y_m})} \right\|_{B\left( {\cc,X} \right)}} \le \frac{C}{{{{(\sum\limits_{j = 1}^m | x - {y_j}|)}^{mn}}}};$
	\item[\emph{(ii)}]
	${\left\| {Q(x,{y_1}, \ldots ,{y_i}, \ldots ,{y_m}) - Q(x,{y_1}, \ldots ,{y_{i'}}, \ldots ,{y_m})} \right\|_{B\left( {\cc,X} \right)}} \le \frac{{C|{y_i} - {y_{i'}}{|^\varepsilon }}}{{{{(\sum\limits_{j = 1}^m | x - {y_j}|)}^{mn + \varepsilon }}}}$\\
	whenever $|x-x'|\leq\frac{1}{2}\max_{1\leq j \leq m}|x-y_j|$;
	\item[\emph{(iii)}]
	${\left\| {Q(x,{y_1}, \ldots ,{y_m}) - Q(x',{y_1}, \ldots ,{y_m})} \right\|_{B\left( {\cc,X} \right)}} \le \frac{{C|x - x'{|^\varepsilon }}}{{{{(\sum\limits_{j = 1}^m | x - {y_j}|)}^{mn + \varepsilon }}}}$\\
	whenever
	$|x-x'|\leq\frac{1}{2}\sum_{j=1}^m|x-y_j|$.
\end{enumerate}
If $\mathop {\min }\limits_{1 \le i \le m} \{ {p_i}\}  > 1$, $\vec{\omega}=(\omega_1,\ldots,\omega_m)\in {A_{\vec P}} \cap {\left( {{A_\infty }}\right)^m}$, then there exists a constant ${C}$, independent of $\vec f$, such that 
\begin{equation*}
{\left\| {T(\vec f)} \right\|_{{L^{p,\kappa}}({v_{\vec \omega }})}} \le C{\prod\limits_{i = 1}^m {\left\| {{f_i}} \right\|} _{{L^{{p_i},\kappa }}({\omega _i})}}.
\end{equation*}
\end{cor}

In theorem \ref{the1}, if we take $X = {L^2}\left( {\left( {0,\infty } \right),\frac{{dt}}{t}} \right)$, we have three important corollaries as follows.

\begin{cor}\label{cor1}
Let $m\in \n$, suppose that $K \in L_{loc}^1(\mm\backslash E),$ $E = \{ (x,\vec y) \in \mm:x = {y_1} =  \cdots  = {y_m}\},$ which satisfies Size Condition:
\begin{equation}\label{N3}
{(\int_0^\infty  {{{\left| {{K_t}(x,\vec y)} \right|}^2}} \frac{{dt}}{t})^{\frac{1}{2}}} \le \frac{C}{{{{(\sum\limits_{j = 1}^m {\left| {x - {y_j}} \right|} )}^{mn}}}}.
\end{equation}
The multilinear square operator $T$ is defined by
\begin{equation}\label{EQQ015}
T(\vec f)(x) = {({\int_0^\infty  {\left| {\int_{\nm} {{K_t}(x,\vec y)\prod\limits_{j = 1}^m {{f_j}({y_j})d{y_1} \cdots d{y_m}} } } \right|} ^2}\frac{{dt}}{t})^{\frac{1}{2}}}.
\end{equation}
The maximal multilinear square operator $T^*$ is defined by
\begin{equation}\label{NW4}
{T^*}(\vec f)(x) = \mathop {\sup }\limits_{\delta  > 0} {T_\delta }(\vec f)(x),
\end{equation}
in which
\begin{equation}
{T_\delta }(\vec f)(x) = {(\int_0^\infty  {{{\left| {\int_{\sum\limits_{j = 1}^m {{{\left| {x - {y_j}} \right|}^2} > {\delta ^2}} } {{K_t}(x,\vec y)\prod\limits_{j = 1}^m {{f_j}({y_j})d{y_1} \cdots d{y_m}} } } \right|}^2}} \frac{{dt}}{t})^{\frac{1}{2}}}.
\end{equation}
Then, under the same operator boundedness condition,  the conclusions in Theorem \ref{the1} still hold for $T$ and $T^*$.
\end{cor}
%
%cor2
%

\begin{cor}\label{cor2}
Let $m\in \n$ and $T$ be an $m$-linear square operator with kernel of C-Z type I or with kernel of type $w(t),$ and $w \in Dini(1)$. If $p_1,\ldots,p_m\in[1,\infty)$ with $1/p=\sum_{k=1}^m 1/{p_k}$, and $\vec{\omega}=(\omega_1,\ldots,\omega_m)\in {A_{\vec P}} \cap {\left( {{A_\infty }}\right)^m}$. For any $0<\kappa<1$, the following results hold:
\begin{enumerate}[(i)]
\item If $\mathop {\min }\limits_{1 \le i \le m} \{ {p_i}\}  > 1,$ then there exists a constant ${C}$, independent of $\vec f$, such that 
\begin{equation}
{\left\| {T(\vec f)} \right\|_{{L^{p,k}}({v_{\vec \omega }})}} \le C{\prod\limits_{i = 1}^m {\left\| {{f_i}} \right\|} _{{L^{{p_i},k}}({\omega _i})}}.
\end{equation}
\item If $\mathop {\min }\limits_{1 \le i \le m} \{ {p_i}\}  = 1,$ then there exists a constant ${C}$, independent of $\vec f$, such that 
\begin{equation}
{\left\| {T(\vec f)} \right\|_{W{L^{p,k}}({v_{\vec \omega }})}} \le C{\prod\limits_{i = 1}^m {\left\| {{f_i}} \right\|} _{{L^{{p_i},k}}({\omega _i})}}.
\end{equation}
\end{enumerate}
In particular, let $K$ be a multilinear Littlewood-Paley (Marcinkiewicz) kernel. Suppose that $0 < \gamma  < \min \{ \delta ,{\gamma _0}\},$ then the above results for multilinear Littlewood-Paley $g$-function is also valid. 
\end{cor}
\begin{cor}\label{cor3}
Let $m\in \n$ and $T_\lambda$ be an $m$-linear square operator with kernel satisfying the integral condition of C-Z type II. If $p_1,\ldots,p_m\in[1,\infty)$ with $1/p=\sum_{k=1}^m 1/{p_k}$, and $\vec{\omega}=(\omega_1,\ldots,\omega_m)\in {A_{\vec P}} \cap {\left( {{A_\infty }} \right)^m}$. For any $0<\kappa<1$ and $\lambda  > 2m$, the following results hold:
\begin{enumerate}[(i)]
\item If $\mathop {\min }\limits_{1 \le i \le m} \{ {p_i}\}  > 1,$ then there exists a constant ${C}$, independent of $\vec f$, such that 
\begin{equation*}
{\left\| {{T_\lambda }(\vec f)} \right\|_{{L^{p,k}}({v_{\vec \omega }})}} \le C{\prod\limits_{i = 1}^m {\left\| {{f_i}} \right\|} _{{L^{{p_i},k}}({\omega _i})}}.
\end{equation*}
\item If $\mathop {\min }\limits_{1 \le i \le m} \{ {p_i}\}  = 1,$ then there exists a constant ${C}$, independent of $\vec f$, such that 
\begin{equation*}
{\left\| {{T_\lambda }(\vec f)} \right\|_{W{L^{p,k}}({v_{\vec \omega }})}} \le C{\prod\limits_{i = 1}^m {\left\| {{f_i}} \right\|} _{{L^{{p_i},k}}({\omega _i})}}.
\end{equation*}
\end{enumerate}
In particular, suppose that $\lambda  > 2m, 0 < \gamma  \le \min \{ \frac{{\lambda n - 2mn}}{2},{\gamma _0},\frac{n}{2}\}$, then the above results for multilinear Littlewood-Paley $g_\lambda ^*$-function are also valid. 
\end{cor}
We also give the similar result for the classical Littlewood-Paley $g$-function on weighted Morrey spaces.
\begin{cor}\label{cor4}
Suppose that $\varphi  \in {L^1}({\rn})$ satisfies 
\begin{enumerate}[(i)]
\item Size Condition:$$\left| {\varphi (x)} \right| \le \frac{C}{{{{(1 + \left| x \right|)}^{n + \alpha }}}}$$
\item Smoothness Condition:$$\left| {\nabla \varphi (x)} \right| \le \frac{C}{{{{(1 + \left| x \right|)}^{n + \alpha '}}}}$$
\item Vanishing Condition:$$\int_{\rn} {\varphi (x)dx = 0} $$ 
\end{enumerate}
for some $\alpha \ge \frac{1}{2}$, $\alpha ' > 1.$
For $1 < p < \infty$, $0<\kappa<1$, if $\omega  \in {A_p}$, then there exists a constant ${C}$, independent of $f$, such that $${\left\| {{g_\varphi }(f)} \right\|_{{L^{p,\kappa}}(\omega )}} \le {C}{\left\| f\right\|_{{L^{p,\kappa}}(\omega )}}.$$
\end{cor}
Next, we consider the boundedness of commutators of the multilinear vector-valued operators on weighted Morrey spaces.

%多线性向量值函数的交换子
\begin{thm}\label{the2}
Let $m\in\n,$ $X$ is a Banach space, $B(\cc,X)$ is the space of all bounded linear operators from $\cc$ to $X$, suppose that operator-valued function $Q:(\mm \backslash E) \to B(\cc,X),$ $E = \{ (x,\vec y) \in \mm:x = {y_1} =  \cdots  = {y_m}\},$ which satisfies Size Condition:
\begin{equation}
	{\left\| {Q(x,\vec y)} \right\|_{B(\cc,X)}} \le \frac{C}{{{{(\sum\limits_{j = 1}^m {\left| {x - {y_j}} \right|} )}^{mn}}}}.
\end{equation}
Set $\vec b = ({b_1}, \cdots ,{b_m}) \in {({BMO})^m},$
we define the commutator of $\vec b$ and multilinear X-valued operator $T$ by
\begin{equation}
	{T_{\vec b}}(\vec f)(x) = \sum\limits_{N = 1}^m {T_{\vec b}^N(\vec f)(x)}.
\end{equation}
where the N-th commutator of $\vec b$ and multilinear X-valued operator $T$ is defined by
\begin{equation}
T_{\vec b}^N(\vec f)(x) = {\left\| {\int_{\nm} {(Q(x,\vec y))(({b_N}(x) - {b_N}({y_N}))\prod\limits_{j = 1}^m {{f_j}({y_j})} )d{y_1} \cdots d{y_m}} } \right\|_X}.
\end{equation}
The iterated commutator of $\vec b$ and multilinear X-valued operator $T$ is defined by
\begin{equation}
{T_{\prod {\vec b} }}(\vec f)(x) = {\left\| {\int_{{\rm{\nm}}} {(Q(x,\vec y))(\prod\limits_{j = 1}^m {({b_j}(x) - {b_j}({y_j}))} {f_j}({y_j}))d{y_1} \cdots d{y_m}} } \right\|_X}.
\end{equation}

If $p_1,\ldots,p_m\in[1,\infty)$ with $1/p=\sum_{k=1}^m 1/{p_k}$, and $\vec{\omega}=(\omega_1,\ldots,\omega_m)\in {A_{\vec P}} \cap {\left( {{A_\infty }} \right)^m}$. For any $0<\kappa<1$, the following results hold:
\begin{enumerate}[(i)]
\item if $\mathop {\min }\limits_{1 \le i \le m} \{ {p_i}\}  > 1,$ such that $T_{\vec b}$ is well-defined on ${L^{{p_1}}}({\omega _1}) \times  \cdots  \times {L^{{p_m}}}({\omega _m})$, which is also bounded from ${L^{{p_1}}}({\omega _1}) \times  \cdots  \times {L^{{p_m}}}({\omega _m})$ to ${L^p}({v_{\vec \omega }}),$ then there exists a constant ${C}$, independent of $\vec f$, such that 
\begin{equation*}
{\left\| {G(\vec f)} \right\|_{{L^{p,\kappa}}({v_{\vec \omega }})}} \le C{\prod\limits_{i = 1}^m {\left\| {{f_i}} \right\|} _{{L^{{p_i},\kappa }}({\omega _i})}}.
\end{equation*}
where $G$ can take $T_{\vec b}$ or $T_{\prod {\vec b} }$;
\item if $\mathop {\min }\limits_{1 \le i \le m} \{ {p_i}\}  = 1,$ such that $T_{\vec b}$ is well-defined on ${L^{{p_1}}}({\omega _1}) \times  \cdots  \times {L^{{p_m}}}({\omega _m})$, which is also bounded from ${L^{{p_1}}}({\omega _1}) \times  \cdots  \times {L^{{p_m}}}({\omega _m})$ to $W{L^p}({v_{\vec \omega }}),$ then there exists a constant ${C}$, independent of $\vec f$, such that  
\begin{equation*}
{\left\| {G(\vec f)} \right\|_{W{L^{p,\kappa}}({v_{\vec \omega }})}} \le C{\prod\limits_{i = 1}^m {\left\| {{f_i}} \right\|} _{{L^{{p_i},\kappa}}({\omega _i})}}.
\end{equation*}
where $G$ can take $T_{\vec b}$ or $T_{\prod {\vec b} }$.
\end{enumerate}
\end{thm}
In theorem \ref{the2}, if we take $X = {L^2}\left( {\left( {0,\infty } \right),\frac{{dt}}{t}} \right)$, we also have a meaningful corollary as follows.
\begin{cor}\label{cor6}
Let $m\in \n$, suppose that $K \in L_{loc}^1(\mm\backslash E),$ $E = \{ (x,\vec y) \in \mm:x = {y_1} =  \cdots  = {y_m}\},$ which satisfies Size Condition:
\begin{equation*}
{(\int_0^\infty  {{{\left| {{K_t}(x,\vec y)} \right|}^2}} \frac{{dt}}{t})^{\frac{1}{2}}} \le \frac{C}{{{{(\sum\limits_{j = 1}^m {\left| {x - {y_j}} \right|} )}^{mn}}}}.
\end{equation*}
Set $\vec b = ({b_1}, \cdots ,{b_m}) \in {({BMO})^m},$
we define the commutator of $\vec b$ and multilinear square operator $T$ by
\begin{equation*}
{T_{\vec b}}(\vec f)(x) = \sum\limits_{N = 1}^m {T_{\vec b}^N(\vec f)(x)},
\end{equation*}
where the N-th commutator of $\vec b$ and multilinear square operator $T$ is defined by
\begin{equation*}
T_{\vec b}^N(\vec f)(x) = {({\int_0^\infty  {\left| {\int_{\nm} {{K_t}(x,\vec y)({b_N}(x) - {b_N}({y_N}))\prod\limits_{j = 1}^m {{f_j}({y_j})} d{y_1} \cdots d{y_m}} } \right|} ^2}\frac{{dt}}{t})^{\frac{1}{2}}}.
\end{equation*}
The iterated commutator of $\vec b$ and multilinear square operator $T$ is defined by
\begin{equation*}
{T_{\prod {\vec b} }}(\vec f)(x) = {(\int_0^\infty  {{{\left| {\int_{{\rm{\nm}}} {{K_t}(x,\vec y)(\prod\limits_{j = 1}^m {({b_j}(x) - {b_j}({y_j}))} {f_j}({y_j}))d{y_1} \cdots d{y_m}} } \right|}^2}} \frac{{dt}}{t})^{\frac{1}{2}}}.
\end{equation*}
Then, under the same operator boundedness condition, the conclusions in Theorem \ref{the2} still hold for $T_{\vec b}$ and $T_{\prod {\vec b} }$.
\end{cor}

According to Corollary \ref{cor6} and Theorem 1.3 in \cite{xue4}, we have a significant corollary for the iterated commutator of multilinear Littlewood–Paley g-function with convolution-type kernel, see \cite{xue4} for more details.
\begin{cor}\label{cor7}
	Let $p_1,\ldots,p_m\in(1,\infty)$ with $1/p=\sum_{k=1}^m 1/{p_k}$, and $\vec{\omega}=(\omega_1,\ldots,\omega_m)\in {A_{\vec P}} \cap {\left( {{A_\infty }} \right)^m}$. Set $\vec b = ({b_1}, \cdots ,{b_m}) \in {({BMO})^m},$ we define the iterated commutator of $\vec b$ and multilinear Littlewood–Paley g-function by
	\begin{equation*}
{{g}_{\prod {\vec b} }}(\vec f)(x) = {(\int_0^\infty  {{{\left| {\int_{{\rm{\nm}}} {{K_t}(x,\vec y)(\prod\limits_{j = 1}^m {({b_j}(x) - {b_j}({y_j}))} {f_j}({y_j}))d{y_1} \cdots d{y_m}} } \right|}^2}} \frac{{dt}}{t})^{\frac{1}{2}}},
	\end{equation*}
	where ${K_t}(x,\vec y) = {\psi _t}\left( {x - {y_1}, \ldots ,x - {y_m}} \right)$.
	Then for any $0<\kappa<1$, there exists a constant ${C}$, independent of $\vec f$, such that  
	\begin{equation*}
	{\left\| {{g_{\prod {\vec b} }}(\vec f)} \right\|_{{L^{p,\kappa }}({v_{\vec \omega }})}} \le C\prod\limits_{i = 1}^m {{{\left\| {{f_i}} \right\|}_{{L^{{p_i},\kappa }}({\omega _i})}}}.
	\end{equation*}
\end{cor}
Similarly, according to Corollary \ref{cor6} and Theorem 1.5 in \cite{He}, we also have a significant corollary for the commutator of multilinear Marcinkiewicz integral with convolution-type homogeneous kernel, see Remark 1.1 in \cite{xue1} and \cite{He}  for more details.
\begin{cor}\label{cor8}
	Let $p_1,\ldots,p_m\in(1,\infty)$ with $1/p=\sum_{k=1}^m 1/{p_k}$, and $\vec{\omega}=(\omega_1,\ldots,\omega_m)\in {A_{\vec P}} \cap {\left( {{A_\infty }} \right)^m}$. Set $\vec b = ({b_1}, \cdots ,{b_m}) \in {({BMO})^m},$ we define the commutator of $\vec b$ and multilinear Marcinkiewicz integral $\mu$ by
	\begin{equation*}
	{\mu_{\vec b}}(\vec f)(x) = \sum\limits_{N = 1}^m {\mu_{\vec b}^N(\vec f)(x)},
	\end{equation*} 
	where the N-th commutator of $\vec b$ and multilinear Marcinkiewicz integral $\mu $ is defined by
	\begin{equation*}
	\mu_{\vec b}^N(\vec f)(x) = {({\int_0^\infty  {\left| {\int_{\nm} {{K_t}(x,\vec y)({b_N}(x) - {b_N}({y_N}))\prod\limits_{j = 1}^m {{f_j}({y_j})} d{y_1} \cdots d{y_m}} } \right|} ^2}\frac{{dt}}{t})^{\frac{1}{2}}},
	\end{equation*}
	where $K(x,\vec y) = \frac{{\Omega \left( {x - {y_1}, \ldots ,x - {y_m}} \right)}}{{{{\left( {\sum\limits_{j = 1}^m {\left| {x - {y_j}} \right|} } \right)}^{m\left( {n - 1} \right)}}}}{\chi _{{{\left( {B\left( {0,1} \right)} \right)}^m}}}\left( {x - {y_1}, \ldots ,x - {y_m}} \right)$.
	Then for any $0<\kappa<1$, there exists a constant ${C}$, independent of $\vec f$ and $\vec b$, such that  
		\begin{equation*}
{\left\| {{\mu _{\vec b}}(\vec f)} \right\|_{{L^{p,\kappa }}({v_{\vec \omega }})}} \le C{\left\| {\vec b} \right\|_{{{(BMO)}^m}}}\prod\limits_{i = 1}^m {{{\left\| {{f_i}} \right\|}_{{L^{{p_i},\kappa }}({\omega _i})}}},
		\end{equation*}
where ${\left\| {{\vec b}} \right\|_{{{(BMO)}^m}}} = \mathop {\sup }\limits_{1 \le N \le m} {\left\| {{b_N}} \right\|_{BMO}}$.
\end{cor}
Now, we give the following weak-type $L\log L$ estimates for the iterated commutator ${T_{\prod {\vec b} }}$ and commutator ${T _{\vec b}}$ of multilinear vector-valued operators on weighted Morrey spaces.
\begin{thm}\label{thm3}
	Let $m\geq2$, $p_i=1$, $i=1,2,\ldots,m$ and $p=1/m$. Set $\vec b = ({b_1}, \cdots ,{b_m}) \in {({BMO})^m}$ and $\vec \omega  \in {A_{(1, \ldots ,1)}} \cap {\left( {{A_\infty }} \right)^m}$, if ${T_{\vec b}^{}}$ and $ {T_{\prod {\vec b} }^{}}$ have weak-type $L\log L$ estimates on weighted Lebesgue spaces, i.e.,
	\begin{equation}\label{eq3}
{\left[ {{v _{\vec \omega}}(\{ x \in {\rn}:|T_{\vec b}^{}(\vec f)(x)| > {\lambda ^m}\} )} \right]^m} \lesssim \Phi ({\left\| {\vec b} \right\|_{{{(BMO)}^m}}})\prod\limits_{k = 1}^m ( \int_{\rn} \Phi  (\frac{{|{f_k}(x)|}}{\lambda }){\omega_k}(x){\mkern 1mu} dx);
	\end{equation}
	\begin{equation}\label{eq4}
{\left[ {{v _{\vec \omega}}(\{ x \in {\rn}:|T_{\prod {\vec b} }^{}(\vec f)(x)| > {\lambda ^m}\} )} \right]^m} \lesssim \prod\limits_{k = 1}^m ( \int_{\rn} {{\Phi ^{(m)}}} (\frac{{|{f_k}(x)|}}{\lambda }){\omega_k}(x){\mkern 1mu} dx),
	\end{equation}
	then, for any given $\lambda>0$ and any ball $B\subset\mathbb R^n$, we have
	\begin{equation*}
	{{{v _{\vec \omega}}{{(B)}^{-m\kappa }}}}{[{{v_{\vec \omega}}(\{ x \in B:|{T _{\vec b}}(\vec f)(x)| > {\lambda ^m}\} )}]^m} \lesssim \Phi ({\left\| {\vec b} \right\|_{{{(BMO)}^m}}})\prod\limits_{i = 1}^m {{{\left\| {\Phi (\frac{{|{f_i}|}}{\lambda })} \right\|}_{{{(L\log L)}^{1,\kappa }}({\omega_i})}}};
	\end{equation*}
	\begin{equation*}
	{{{v _{\vec \omega}}{{(B)}^{-m\kappa }}}}{[{v_{\vec \omega}}(\{ x \in B:\left| {{T_{\prod {\vec b} }}(\vec f)(x)} \right| > {\lambda ^m}\} )]^m}
	\lesssim
	\prod_{i=1}^m\bigg\|\Phi^{(m)}\bigg(\frac{|f_i|}{\lambda}\bigg)\bigg\|_{(L\log L)^{1,\kappa}(\omega_i)},
	\end{equation*}
	where ${\left\| {{\vec b}} \right\|_{{{(BMO)}^m}}} = \mathop {\sup }\limits_{1 \le N \le m} {\left\| {{b_N}} \right\|_{BMO}}$, $\Phi(t):=t(1+\log^+t)$, $\log^+t:=\max\{\log t,0\}$ and $\Phi^{(m)}=\overbrace{\Phi\circ\cdots\circ\Phi}^m$.
\end{thm}

Finally, combining Theorem \ref{thm3}, Theorem 3.16 in \cite{Lerner}, Theorem 1.6 in \cite{He} and Theorem 1.4 in \cite{xue4}, we give the following weak-type $L\log L$ estimates for iterated commutator ${g_{\prod {\vec b} }}$ and commutator ${\mu_{\vec b}}$.
\begin{cor}
	Let $m\geq2$, $p_i=1$, $i=1,2,\ldots,m$ and $p=1/m$. Set $\vec b = ({b_1}, \cdots ,{b_m}) \in {({BMO})^m}$ and $\vec \omega  \in {A_{(1, \ldots ,1)}} \cap {\left( {{A_\infty }} \right)^m}$, then, for any given $\lambda>0$ and any ball $B\subset\mathbb R^n$, we have
	\begin{equation*}
	{{{v_{\vec \omega}}{{(B)}^{-m\kappa }}}}{[{{v_{\vec \omega}}(\{ x \in B:|{\mu _{\vec b}}(\vec f)(x)| > {\lambda ^m}\} )}]^m} \lesssim \Phi ({\left\| {\vec b} \right\|_{{{(BMO)}^m}}})\prod\limits_{i = 1}^m {{{\left\| {\Phi (\frac{{|{f_i}|}}{\lambda })} \right\|}_{{{(L\log L)}^{1,\kappa }}({\omega_i})}}};
	\end{equation*}
	\begin{equation*}
{{{v_{\vec \omega}}{{(B)}^{-m\kappa }}}}{[{v_{\vec \omega}}(\{ x \in B:\left| {{g_{\prod {\vec b} }}(\vec f)(x)} \right| > {\lambda ^m}\} )]^m}
	\lesssim
	\prod_{i=1}^m\bigg\|\Phi^{(m)}\bigg(\frac{|f_i|}{\lambda}\bigg)\bigg\|_{(L\log L)^{1,\kappa}(\omega_i)},
	\end{equation*}
where ${\left\| {{\vec b}} \right\|_{{{(BMO)}^m}}} = \mathop {\sup }\limits_{1 \le N \le m} {\left\| {{b_N}} \right\|_{BMO}}$, $\Phi(t):=t(1+\log^+t)$, $\log^+t:=\max\{\log t,0\}$ and $\Phi^{(m)}=\overbrace{\Phi\circ\cdots\circ\Phi}^m$.
\end{cor}

The organization of this paper is as follows. In section \ref{sec2}, we prepare some definitions and preliminary lemmas, which play a fundamental role in this paper. Section \ref{sec3} is the proofs of our main results. References are given at the end of the paper.

Throughout this article, we will use $C$ to denote a positive constant, which is independent of the main parameters and not necessarily the same at each occurrence. 

By $A\lesssim B$, we mean that there exists a constant $C>0$, such that $A \le CB$.

By $ A \approx B$, we mean that $A\lesssim B$ and $B\lesssim A$.
%------------------------section 2------------------------
\section{Preliminaries}\label{sec2}
\par
First let us recall some standard definitions and notations. The classical $A_p$ weight
theory was introduced by Muckenhoupt in the study of weighted
$L^p$ boundedness of Hardy-Littlewood maximal functions, one can see Chapter 7 in \cite{Gra1}.
\begin{defn}[\cite{Gra1}]
A weight $\omega$ is a nonnegative locally integrable function on $\mathbb R^n$ that takes values in $(0,\infty)$ almost everywhere. We denote the ball with the center $x_0$ and radius $r$ by $B=B(x_0,r)$, we say that $\omega\in A_p$,\,$1<p<\infty$, if
$$\left(\frac1{|B|}\int_B \omega(x)\,dx\right)\left(\frac1{|B|}\int_B \omega(x)^{-\frac{1}{p-1}}\,dx\right)^{p-1}\le C \quad\mbox{for every ball}\; B\subseteq \mathbb
R^n,$$ where $C$ is a positive constant which is independent of $B$.\\
We say $\omega\in A_1$, if
$$\frac1{|B|}\int_B \omega(x)\,dx\le C\,\underset{x\in B}{\mbox{ess\,inf}}\,\omega(x)\quad\mbox{for every ball}\;B\subseteq\mathbb R^n.$$
We denote $${A_\infty } = \bigcup\limits_{1 \le p < \infty } {{A_p}}.$$
\end{defn}
\begin{defn}[\cite{Gra1}]
A weight function $\omega$ is said to belong to the reverse H\"{o}lder class $RH_r$ if there exist two constants $r>1$ and $C>0$ such that the following reverse H\"{o}lder inequality holds
$$\left(\frac{1}{|B|}\int_B \omega(x)^r\,dx\right)^{1/r}\le C\left(\frac{1}{|B|}\int_B \omega(x)\,dx\right)\quad\mbox{for every ball}\; B\subseteq \mathbb R^n.$$
It is well known that if $\omega\in A_p$ with $1<p<\infty$, then $\omega\in A_r$ for all $r>p$, and $\omega\in A_q$ for some $1<q<p$. If $\omega\in A_p$ with $1\le p<\infty$, then there exists $r>1$ such that $\omega\in RH_r$.
\end{defn}
Now let us recall the definitions of multiple weights. 
\begin{defn}[\cite{Lerner}]
For $m$ exponents $p_1,\ldots,p_m$, we will write $\vec{P}$ for the vector $\vec{P}=(p_1,\ldots,p_m)$. Let $p_1,\ldots,p_m\in[1,\infty)$ and $p\in(0,\infty)$ with $1/p=\sum_{k=1}^m 1/{p_k}$. Given $\vec{\omega}=(\omega_1,\ldots,\omega_m)$, set $v_{\vec{\omega}}=\prod_{i=1}^m \omega_i^{p/{p_i}}$. We say that $\vec{\omega}$ satisfies the $A_{\vec{P}}$ condition if it satisfies
\begin{equation}
\sup_B\left(\frac{1}{|B|}\int_B v_{\vec{\omega}}(x)\,dx\right)^{1/p}\prod_{i=1}^m\left(\frac{1}{|B|}\int_B \omega_i(x)^{1-p'_i}\,dx\right)^{1/{p'_i}}<\infty.
\end{equation}
when $p_i=1,$ $\left(\frac{1}{|B|}\int_B \omega_i(x)^{1-p'_i}\,dx\right)^{1/{p'_i}}$ is understood as ${(\mathop {\inf }\limits_{x \in B} {\omega _i}(x))^{ - 1}}$.
\end{defn}

\begin{lem}[\cite{Lerner}]\label{lem4}
	Let $p_1,\ldots,p_m\in[1,+\infty)$ and $1/p=\sum_{k=1}^m 1/{p_k}$. Then $\vec{\omega}=(\omega_1,\ldots,\omega_m)\in A_{\vec{P}}$ if and only if
	\begin{equation}\label{multi2}
	\left\{
	\begin{aligned}
	&v_{\vec{\omega}}\in A_{mp},\\
	&\omega_k^{1-p'_k}\in A_{mp'_k},\quad k=1,\ldots,m,
	\end{aligned}\right.
	\end{equation}
	where $\nu_{\vec{\omega}}=\prod_{k=1}^m \omega_k^{p/{p_k}}$ and the condition $\omega_k^{1-p'_k}\in A_{mp'_k}$ in the case $p_k=1$ is understood as $\omega_k^{1/m}\in A_1$.
\end{lem}

\par
Given a ball $B$ and $\lambda>0$, $\lambda B$ denotes the ball with the same center as $B$ whose radius is $\lambda$ times that of $B$. For a given weight function $\omega$ and a measurable set $E$, we also denote the Lebesgue measure of $E$ by $|E|$ and the weighted measure of $E$ by $\omega(E)$, where $\omega(E)=\int_E \omega(x)\,dx$. 

\begin{defn}[\cite{Wang4}]
Given a Young function $\Phi$ and $\omega\in A_\infty$, we define
\begin{equation*}
\big\|f\big\|_{\Phi(\omega),B}:=\inf\bigg\{\sigma>0:\frac{1}{\omega(B)}
\int_B\Phi\bigg(\frac{|f(x)|}{\sigma}\bigg)\cdot \omega(x)\,dx\leq1\bigg\}.
\end{equation*}
When $\Phi(t)=t$, this norm is denoted by $\|\cdot\|_{L(\omega),B}$, when $\Phi(t)=t(1+\log^+t)$, this norm is denoted by $\|\cdot\|_{L\log L(\omega),B}$. The complementary Young function of $\Phi(t)$ is $\bar{\Phi}(t)\approx\exp(t)-1$ with the norm denoted by $\|\cdot\|_{\exp L(\omega),B}$. For $\omega\in A_\infty$ and any $B$ in $\mathbb R^n$, the following generalized H\"older's inequality is valid.
\begin{equation}\label{Wholder}
\frac{1}{\omega(B)}\int_B\big|f(x)\cdot g(x)\big|\omega(x)\,dx\leq C\big\|f\big\|_{L\log L(\omega),B}\big\|g\big\|_{\exp L(\omega),B}.
\end{equation}
\end{defn}
	
Now, we give the following results that we will use frequently in the sequel.
\begin{lem}[\cite{Gra1}]\label{lem1}
Let $\omega\in A_p$, $p\ge1$. Then, for any ball $B$, there exists an absolute constant $C$ such that
\begin{equation*}
\omega(2B)\le C \omega(B).
\end{equation*}
In general, for any $\lambda>1$, we have
\begin{equation}\label{LEQ1}
\omega(\lambda B)\le C\lambda^{np}\omega(B)
\end{equation}
where $C$ does not depend on $B$ nor on $\lambda$.
\end{lem}

\begin{lem}[\cite{Gra2}]\label{lem2}
For all $p \in [1,\infty )$ and $f \in L_{loc}^1(\rn)$, we have
\begin{equation}
\mathop {\sup }\limits_B {(\frac{1}{{\left| B \right|}}\int_B^{} {{{\left| {f(x) - {f_B}} \right|}^p}dx} )^{\frac{1}{p}}} \approx {\left\| f \right\|_{BMO}}:=\mathop {\sup }\limits_B (\frac{1}{{\left| B \right|}}\int_B^{} {\left| {f(x) - {f_B}} \right|dx} ).
\end{equation}
\end{lem}

\begin{lem}[\cite{Gra2}]
For all $p \in (0,\infty )$ and $f \in BMO$, we have
\begin{equation}\label{eq6}
\mathop {\sup }\limits_B {(\frac{1}{{\left| B \right|}}\int_B^{} {{{\left| {f(x) - {f_B}} \right|}^p}dx} )^{\frac{1}{p}}} \lesssim {\left\| f \right\|_{BMO}}.
\end{equation}
\end{lem}

\begin{lem}[\cite{Wang4}]\label{BMO3}
	Let $b$ be a function in $\mathrm{BMO}(\mathbb R^n)$. Then for any ball $B$ in $\mathbb R^n$ and any $\omega\in A_{\infty}$, we have
	\begin{equation}\label{BMOwang}
	\big\|b-b_{B}\big\|_{\exp L(\omega),B}\lesssim \|b\|_{BMO}.
	\end{equation}
\end{lem}

\begin{lem}[\cite{Gra1}]\label{lem3}
Let $\omega\in RH_r$ with $r>1$. Then there exists a constant $C$ such that
\begin{equation}\label{AAA5}
\frac{\omega(E)}{\omega(B)}\le C\left(\frac{|E|}{|B|}\right)^{(r-1)/r}
\end{equation}
for any measurable subset $E$ of a ball $B$.
\end{lem}
\begin{lem}\label{lem5}
	Let $m\in \n,$ $p_1,\ldots,p_m\in[1,\infty)$ with $1/p=\sum_{k=1}^m 1/{p_k}$. Assume that $\omega_1,\ldots,\omega_m\in A_\infty$ and $v_{\vec{\omega}}=\prod_{i=1}^m \omega_i^{p/{p_i}}\in A_\infty$, then for any ball $B,$ we have
	\begin{equation*}\label{N1}
	\prod\limits_{i = 1}^m {{{\left( {\int_B {{\omega _i}} (x){\mkern 1mu} dx} \right)}^{p/{p_i}}}}  \approx \int_B {{v _{\vec \omega }}} (x){\mkern 1mu} dx.
	\end{equation*}
	\begin{proof}[Proof:]
	Using Jensen's inequality and the definition of ${A_\infty }$ which can be found in \cite[p. 12]{Gra1} and \cite[p. 525]{Gra1}, we get
	\begin{equation*}
	\left| B \right|\exp (\frac{1}{{\left| B \right|}}\int_B {\log \omega _i^{{q_i}}} ) \le \omega _i^{{q_i}}(B) \lesssim \left| B \right|\exp (\frac{1}{{\left| B \right|}}\int_B {\log \omega _i^{{q_i}}} ).
	\end{equation*}
	and then we have
	\begin{align*}
	\prod\limits_{i = 1}^m {\omega _i^{{q_i}}{{(B)}^{\frac{q}{{{q_i}}}}}}  \approx \left| B \right|\exp (\frac{1}{{\left| B \right|}}\int_B {\log u_{\vec \omega }^q} ) \approx u_{\vec \omega }^q(B).
	\end{align*}
\end{proof}
\end{lem}
Given a weight function $\omega$ on $\mathbb R^n$, for $0<p<\infty$, we denote by $f \in {L^p}(X,\omega dx)$ the space of all functions satisfying
$${\left\| f \right\|_{{L^p}(X,\omega dx)}}: = {(\int_X^{} {{{\left| {f(x)} \right|}^p}} \omega (x)dx)^{\frac{1}{p}}} < \infty.$$
For simplicity, we abbreviate ${L^p}(\rn,\omega dx)$ to ${L^p}(\omega )$.\\
For $0<p<\infty$, we also denote by $W{L^p}(X,\omega dx)$ the weighted weak Lebesgue space consisting of all measurable functions $f$ satisfying
\begin{equation*}
{\left\| f \right\|_{W{L^p}(X,\omega dx)}} = \mathop {\sup }\limits_{\lambda  > 0} \lambda  \cdot \omega {(\{ x \in X:|f(x)| > \lambda \} )^{1/p}} < \infty.
\end{equation*}
For simplicity, we abbreviate $W{L^p}(\rn,\omega dx)$ to $W{L^p}(\omega )$.\\
\par Before proving the main theorems, we give some useful reults as follows.
\begin{lem}[\cite{LuDingYan}]\label{lem6}
Suppose that $\varphi  \in {L^1}({\rn})$ satisfies 
\begin{enumerate}[(i)]
\item Size Condition:$$\left| {\varphi (x)} \right| \le \frac{B_1}{{{{(1 + \left| x \right|)}^{n + \alpha }}}}$$
\item Smoothness Condition:$$\left| {\nabla \varphi (x)} \right| \le \frac{B_2}{{{{(1 + \left| x \right|)}^{n + \alpha '}}}}$$
\item Vanishing Condition:$$\int_{\rn} {\varphi (x)dx = 0} $$ 
\end{enumerate}
for some $\alpha  > 0$, $\alpha ' > 1.$
For $1 < p < \infty$, if $\omega  \in {A_p}$, then we have:
\begin{equation}\label{EQ004}
{\left\| {{g_\varphi }(f)} \right\|_{{L^p}(\omega )}} \mathbin{\lower.3ex\hbox{$\buildrel<\over
{\smash{\scriptstyle\sim}\vphantom{_x}}$}} {\left\| f \right\|_{{L^p}(\omega )}}
\end{equation}
\end{lem}
\begin{lem}[\cite{xue1}]\label{lem7}
Let $m\in \n$ and $T$ be an $m$-linear square operator with kernel satisfying the integral condition of C-Z type I or type $w$, and $w \in Dini(1)$. If $p_1,\ldots,p_m\in[1,\infty)$, and $p\in(0,\infty)$ with $1/p=\sum_{k=1}^m 1/{p_k}$, and $\vec{\omega}=(\omega_1,\ldots,\omega_m)\in A_{\vec{P}}$, the following results hold:
\begin{enumerate}[(i)]
\item If $\mathop {\min }\limits_{1 \le i \le m} \{ {p_i}\}  > 1,$ then there exists a constant ${C}$, independent of $\vec f$, such that 
\begin{equation}
{\left\| {T(\vec f)} \right\|_{{L^{p}}({v_{\vec \omega }})}} \le C{\prod\limits_{i = 1}^m {\left\| {{f_i}} \right\|} _{{L^{{p_i}}}({\omega _i})}}.
\end{equation}
\item If $\mathop {\min }\limits_{1 \le i \le m} \{ {p_i}\}  = 1,$ then there exists a constant ${C}$, independent of $\vec f$, such that 
\begin{equation}
{\left\| {T(\vec f)} \right\|_{W{L^{p}}({v_{\vec \omega }})}} \le C{\prod\limits_{i = 1}^m {\left\| {{f_i}} \right\|} _{{L^{{p_i}}}({\omega _i})}}.
\end{equation}
\end{enumerate}
\end{lem}

\begin{lem}[\cite{xue1}]\label{lem8}
Let $m\in \n$ and $T_\lambda$ be an $m$-linear square operator with kernel satisfying the integral condition of C-Z type II. If $p_1,\ldots,p_m\in[1,\infty)$, and $p\in(0,\infty)$ with $1/p=\sum_{k=1}^m 1/{p_k}$, and $\vec{\omega}=(\omega_1,\ldots,\omega_m)\in A_{\vec{P}}$, then for any $\lambda  > 2m$, the following results hold:
\begin{enumerate}[(i)]
\item If $\mathop {\min }\limits_{1 \le i \le m} \{ {p_i}\}  > 1,$ then there exists a constant ${C}$, independent of $\vec f$, such that 
\begin{equation}\label{LEQ2}
{\left\| {{T_\lambda }(\vec f)} \right\|_{{L^{p}}({v_{\vec \omega }})}} \le C{\prod\limits_{i = 1}^m {\left\| {{f_i}} \right\|} _{{L^{{p_i}}}({\omega _i})}}.
\end{equation}
\item If $\mathop {\min }\limits_{1 \le i \le m} \{ {p_i}\}  = 1,$ then there exists a constant ${C}$, independent of $\vec f$, such that 
\begin{equation}\label{LEQ3}
{\left\| {{T_\lambda }(\vec f)} \right\|_{W{L^{p}}({v_{\vec \omega }})}} \le C{\prod\limits_{i = 1}^m {\left\| {{f_i}} \right\|} _{{L^{{p_i}}}({\omega _i})}}.
\end{equation}
\end{enumerate}
\end{lem}
\vspace{0.1cm}

\vspace{0.1cm}

%------------------------section 3------------------------
\section{Proofs of Main Results}\label{sec3}

\subsection{Proof of Theorem \ref{the1}}
\subsubsection{Proof of (i) of Theorem \ref{the1}}
\begin{proof}[Proof:]
For any ball $B=B(x_0,r)$, let $f_i=f^0_i+f^{\infty}_i$, where $f^0_i=f_i\chi_{2B}$, $i=1,\ldots,m$ and $\chi_{2B}$ denotes the characteristic function of $2B$. Then we write
\begin{equation*}
\begin{split}
\prod_{i=1}^m f_i(y_i)&=\prod_{i=1}^m\Big(f^0_i(y_i)+f^{\infty}_i(y_i)\Big)\\
&=\sum_{\alpha_1,\ldots,\alpha_m\in\{0,\infty\}}f^{\alpha_1}_1(y_1)\cdots f^{\alpha_m}_m(y_m)\\
&=\prod\limits_{i = 1}^m {f_i^0} ({y_i}) + \sum\limits_{{\alpha _1} +  \cdots  + {\alpha _m} \ne 0} {f_1^{{\alpha _1}}({y_1}) \cdots f_m^{{\alpha _m}}({y_m})}.
\end{split}
\end{equation*}
Since $T$ is an $m$-sublinear operator, then we have
\begin{equation*}
\begin{split}
&{\nu _{\vec \omega }}{(B)^{ - \frac{\kappa }{p}}}{\left\| {T({f_1}, \cdots ,{f_m})} \right\|_{{L^p}(B,{\nu _{\vec \omega }}dx)}}\\
\le &{\nu _{\vec \omega }}{(B)^{ - \frac{\kappa }{p}}}{\left\| {T(f_1^0, \cdots ,f_m^0)} \right\|_{{L^p}(B,{\nu _{\vec \omega }}dx)}} + \sum\limits_{{\alpha _1} +  \cdots  + {\alpha _m} \ne 0}^{} {{\nu _{\vec \omega }}{{(B)}^{ - \frac{\kappa }{p}}}{{\left\| {T(f_1^{{\alpha _1}}, \cdots ,f_m^{{\alpha _m}})} \right\|}_{_{{L^p}(B,{\nu _{\vec \omega }}dx)}}}} \\
: = &{I^0} + \sum\limits_{{\alpha _1} +  \cdots  + {\alpha _m} \ne 0}^{} {{I^{{\alpha _1} \cdots ,{\alpha _m}}}}.
\end{split}
\end{equation*}
For finishing the proof, we merely need to prove:
\begin{equation}\label{N2}
{I^{{\alpha _1}, \cdots ,{\alpha _m}}} \lesssim \prod\limits_{i = 1}^m {{{\left\| {{f_i}} \right\|}_{{L^{{p_i},\kappa }}({\omega _i})}}}, 
\end{equation}
where ${\alpha _i} \in \{ 0,\infty \} ,i = 1, \cdots ,m.$\\
In view of Lemma \ref{lem4}, we have $\nu_{\vec{\omega}}\in A_{mp}$. Applying the boundedness, Lemma \ref{lem1} and Lemma \ref{lem5}, we get
\begin{equation*}
\begin{split}
I^0&\lesssim
\frac{1}{\nu_{\vec{\omega}}(B)^{\kappa/p}}\prod_{i=1}^m\left(\int_{2B}|f_i(x)|^{p_i}\omega_i(x)\,dx\right)^{1/{p_i}}\\
&\lesssim\prod_{i=1}^m\big\|f_i\big\|_{L^{p_i,\kappa}(\omega_i)}\cdot
\frac{\prod_{i=1}^m \omega_i(2B)^{\kappa/{p_i}}}{\nu_{\vec{\omega}}(B)^{\kappa/p}}\\
&\lesssim\prod_{i=1}^m\big\|f_i\big\|_{L^{p_i,\kappa}(\omega_i)}\cdot
\frac{\nu_{\vec{\omega}}(2B)^{\kappa/p}}{\nu_{\vec{\omega}}(B)^{\kappa/p}}\\
&\lesssim\prod_{i=1}^m \big\|f_i\big\|_{L^{p_i,\kappa}(\omega_i)}.
\end{split}
\end{equation*}

For $T$, we have 
\begin{equation*}
\begin{split}
T(\vec f)(x) &= {\left\| {\int_{\nm} {(Q(x,\vec y))(\prod\limits_{j = 1}^m {{f_j}({y_j}))d{y_1} \cdots d{y_m}} } } \right\|_X}\\
&\le \int_{\nm} {{{\left\| {(Q(x,\vec y))(\prod\limits_{j = 1}^m {{f_j}({y_j})} )} \right\|}_X}d{y_1} \cdots d{y_m}}\\
&\lesssim \int_{\nm} {\frac{{\left| {\prod\limits_{j = 1}^m {{f_j}({y_j})} } \right|}}{{{{(\sum\limits_{i = 1}^m {\left| {x - {y_i}} \right|} )}^{mn}}}}d{y_1} \cdots d{y_m}}.
\end{split}
\end{equation*}

To obtain the conclusions, we establish some geometric relationships by trigonometric inequality as follows:
\begin{enumerate}[(i)]
\item If $x\in B, y\in (2B)^c$, we obviously have:$$\left| {x - y} \right| \approx \left| {{x_0} - y} \right|;$$ 
\item If $x\in B$, $y \in {2^{j + 1}}B\backslash {2^j}B$, $j \in \n,$ we obviously have:$${2^{j - 1}}r \le \left| {x - y} \right| \le {2^{j + 2}}r.$$
\end{enumerate}
For the other terms, we first consider the case when $\alpha_1=\cdots=\alpha_m=\infty$. For $x\in B,$ we have
\begin{align}
\big|T(f^\infty_1,\ldots,f^\infty_m)(x)\big|
&\lesssim \int_{(\mathbb R^n)^m\backslash(2B)^m}
\frac{|f_1(y_1)\cdots f_m(y_m)|}{(|x-y_1|+\cdots+|x-y_m|)^{mn}}dy_1\cdots dy_m\notag\\
&= \sum_{j=1}^\infty\int_{(2^{j+1}B)^m\backslash(2^{j}B)^m}
\frac{|f_1(y_1)\cdots f_m(y_m)|}{(|x-y_1|+\cdots+|x-y_m|)^{mn}}dy_1\cdots dy_m\notag\\
&\lesssim \sum\limits_{j = 1}^\infty  {\int_{{{({2^{j + 1}}B)}^m}\backslash {{({2^j}B)}^m}}^{} {\frac{{|{f_1}({y_1}) \cdots {f_m}({y_m})|}}{{{{({2^{j + 1}}r)}^{mn}}}}d{y_1} \cdots d{y_m}} }\notag\\
&\lesssim \sum_{j=1}^\infty\prod_{i=1}^m\frac{1}{|2^{j+1}B|}\int_{2^{j+1}B}\big|f_i(y_i)\big|\,dy_i.\label{AA1}
\end{align}
By using H\"older's inequality, the $A_{\vec{P}}$ condition and Lemma \ref{lem5}, we obtain:
\begin{equation*}
\begin{split}
\big|T(f^\infty_1,\ldots,f^\infty_m)(x)\big|&\lesssim \sum_{j=1}^\infty\prod_{i=1}^m
\frac{1}{|2^{j+1}B|}\left(\int_{2^{j+1}B}\big|f_i(y_i)\big|^{p_i}\omega_i(y_i)\,dy_i\right)^{1/{p_i}}
\left(\int_{2^{j+1}B}\omega_i(y_i)^{1-p'_i}\,dy_i\right)^{1/{p'_i}}\\
&\lesssim\sum_{j=1}^\infty\frac{1}{|2^{j+1}B|^m}\cdot
\frac{|2^{j+1}B|^{\frac 1p+\sum_{i=1}^m(1-\frac 1{p_i})}}{\nu_{\vec{\omega}}(2^{j+1}B)^{1/p}}
\prod_{i=1}^m
\left(\big\|f_i\big\|_{L^{p_i,\kappa}(\omega_i)}\omega_i\big(2^{j+1}B\big)^{\kappa/{p_i}}\right)\\
&=(\prod_{i=1}^m \big\|f_i\big\|_{L^{p_i,\kappa}(\omega_i)})\cdot\sum_{j=1}^\infty
\left(\frac{\prod_{i=1}^m \omega_i(2^{j+1}B)^{\kappa/{p_i}}}{\nu_{\vec{\omega}}(2^{j+1}B)^{1/p}}\right)\\
&\lesssim(\prod_{i=1}^m \big\|f_i\big\|_{L^{p_i,\kappa}(\omega_i)})\cdot
\sum_{j=1}^\infty\nu_{\vec{\omega}}\big(2^{j+1}B\big)^{{(\kappa-1)}/p}.
\end{split}
\end{equation*}
Thus, we have
\begin{equation*}
\begin{split}
I^{\infty,\ldots,\infty}&\lesssim(\prod_{i=1}^m \big\|f_i\big\|_{L^{p_i,\kappa}(\omega_i)})\cdot\sum_{j=1}^\infty
\frac{\nu_{\vec{\omega}}(B)^{{(1-\kappa)}/p}}{\nu_{\vec{\omega}}(2^{j+1}B)^{{(1-\kappa)}/p}}\\
&\lesssim(\prod_{i=1}^m \big\|f_i\big\|_{L^{p_i,\kappa}(\omega_i)})\cdot\sum_{j=1}^\infty
\left(\frac{|B|}{|2^{j+1}B|}\right)^{{\delta(1-\kappa)}/p}\\
&\lesssim\prod_{i=1}^m \big\|f_i\big\|_{L^{p_i,\kappa}(\omega_i)},
\end{split}
\end{equation*}
where we use the fact as follow:
\begin{equation*}\label{EQQ009} 
\frac{\nu_{\vec{w}}(B)}{\nu_{\vec{w}}(2^{j+1}B)}\lesssim\left(\frac{|B|}{|2^{j+1}B|}\right)^\delta,
\end{equation*}
since we know $\nu_{\vec{w}}\in A_{mp} \subseteq A_\infty,$ and apply Lemma \ref{lem3}. The last inequality holds since $0<\kappa<1$ and $\delta>0$.\\
Without loss of generality, we may assume that ${\alpha _1} =  \cdots  = {\alpha _\ell } = \infty $, and ${\alpha _{l + 1}} =  \cdots  = {\alpha _m} = 0$. For any $x\in B$, we have 
\begin{align}
&\big|T(f^\infty_1,\ldots,f^\infty_\ell,f^0_{\ell+1},\ldots,f^0_m)(x)\big|\notag\\
\lesssim&\int_{(\mathbb R^n)^{\ell}\backslash(2B)^{\ell}}\int_{(2B)^{m-\ell}}\frac{|f_1(y_1)\cdots f_m(y_m)|}{(|x-y_1|+\cdots+|x-y_m|)^{mn}}dy_1\cdots dy_m\notag\\
\lesssim&(\prod_{i=\ell+1}^m\int_{2B}\big|f_i(y_i)\big|\,dy_i)\times\sum_{j=1}^\infty\frac{1}{|2^{j+1}B|^m}\int_{(2^{j+1}B)^\ell\backslash(2^{j}B)^\ell}
\big|f_1(y_1)\cdots f_{\ell}(y_\ell)\big|\,dy_1\cdots dy_\ell\notag\\
\le& \sum_{j=1}^\infty\prod_{i=1}^m\frac{1}{|2^{j+1}B|}\int_{2^{j+1}B}\big|f_i(y_i)\big|\,dy_i,\label{EQQ011}
\end{align}
where the second inequality is valid, since the calculations here are similar to before.
It is the same situation as before, so for any $x\in B$, we also have
\begin{equation}\label{EQQ012}
\big|T(f^\infty_1,\ldots,f^\infty_\ell,f^0_{\ell+1},\ldots,f^0_m)(x)\big|\lesssim(\prod_{i=1}^m \big\|f_i\big\|_{L^{p_i,\kappa}(\omega_i)})\cdot
\sum_{j=1}^\infty\nu_{\vec{\omega}}\big(2^{j+1}B\big)^{{(\kappa-1)}/p}.
\end{equation}
Consequently, we finish the proof by 
\begin{equation*}
\begin{split}
{I^{\infty , \ldots \infty ,0, \ldots ,0}}&\le\nu_{\vec{\omega}}(B)^{{(1-\kappa)}/p}
\big|T(f^\infty_1,\ldots,f^\infty_\ell,f^0_{\ell+1},\ldots,f^0_m)(x)\big|\\
&\lesssim(\prod_{i=1}^m \big\|f_i\big\|_{L^{p_i,\kappa}(\omega_i)})\cdot\sum_{j=1}^\infty
\frac{\nu_{\vec{\omega}}(B)^{{(1-\kappa)}/p}}{\nu_{\vec{\omega}}(2^{j+1}B)^{{(1-\kappa)}/p}}\\
&\lesssim(\prod_{i=1}^m \big\|f_i\big\|_{L^{p_i,\kappa}(\omega_i)})\cdot\sum_{j=1}^\infty
\left(\frac{|B|}{|2^{j+1}B|}\right)^{{\delta(1-\kappa)}/p}\\
&\lesssim\prod_{i=1}^m \big\|f_i\big\|_{L^{p_i,\kappa}(\omega_i)}.
\end{split}
\end{equation*}
Combining with (\ref{N2}), we have already finished this proof.
\end{proof}

\subsubsection{Proof of (ii) of Theorem \ref{the1}}
\begin{proof}[Proof:]
For any ball $B=B(x_0,r)$, decompose $f_i=f^0_i+f^{\infty}_i$, where $f^0_i=f_i\chi_{2B}$, $i=1,\ldots,m$. For each $\lambda>0$, we have
\begin{equation*}
\begin{split}
&{\nu _{\vec \omega }}{(B)^{ - \frac{\kappa }{p}}}{\left\| {T({f_1}, \cdots ,{f_m})} \right\|_{W{L^p}(B,{\nu _{\vec \omega }}dx)}}\\
\lesssim&{\nu _{\vec \omega }}{(B)^{ - \frac{\kappa }{p}}}{\left\| {T(f_1^0, \cdots ,f_m^0)} \right\|_{W{L^p}(B,{\nu _{\vec \omega }}dx)}} + \sum\limits_{{\alpha _1} +  \cdots  + {\alpha _m} \ne 0}^{} {{\nu _{\vec \omega }}{{(B)}^{ - \frac{\kappa }{p}}}{{\left\| {T(f_1^{{\alpha _1}}, \cdots ,f_m^{{\alpha _m}})} \right\|}_{_{W{L^p}(B,{\nu _{\vec \omega }}dx)}}}}\\
: = &{J^0} + \sum\limits_{{\alpha _1} +  \cdots  + {\alpha _m} \ne 0} {{J^{{\alpha _1} \cdots ,{\alpha _m}}}}.
\end{split}
\end{equation*}
For finishing the proof, we merely need to prove:
\begin{equation}\label{N4}
{J^{{\alpha _1}, \cdots ,{\alpha _m}}} \lesssim \prod_{i=1}^m \big\|f_i\big\|_{L^{p_i,\kappa}(\omega_i)}.
\end{equation}
where ${\alpha _i} \in \{ 0,\infty \} ,i = 1, \cdots ,m.$\\
We know that $\nu_{\vec{\omega}}\in A_{mp}$ with $1\le mp<\infty$.
Similarly, we have
\begin{equation*}
\begin{split}
J^0&\lesssim{\nu _{\vec \omega }}{(B)^{ - \frac{\kappa }{p}}}{\prod\limits_{i = 1}^m {\left\| {{f_i}} \right\|} _{{L^{{p_i}}}(2B,{\omega _i}dx)}}\\
&\lesssim\frac{{\prod\limits_{i = 1}^m {{\omega_i}} {{(2B)}^{\kappa /{p_i}}}}}{{{\nu _{\vec \omega }}{{(B)}^{\kappa /p}}}}{\prod\limits_{i = 1}^m {\left\| {{f_i}} \right\|} _{{L^{{p_i},\kappa }}({\omega _i})}}\\
&\lesssim\frac{{{\nu _{\vec \omega}}{{(2B)}^{\kappa /p}}}}{{{\nu _{\vec \omega }}{{(B)}^{\kappa /p}}}}{\prod\limits_{i = 1}^m {\left\| {{f_i}} \right\|} _{{L^{{p_i},\kappa }}({\omega _i})}}\\
&\lesssim\prod_{i=1}^m \big\|f_i\big\|_{L^{p_i,\kappa}(\omega_i)}.
\end{split}
\end{equation*}
Consequently, we have ${J^0} \lesssim \prod_{i=1}^m \big\|f_i\big\|_{L^{p_i,\kappa}(\omega_i)}.$\\
In proof of (i), for ${\alpha _1} +  \cdots  + {\alpha _m} \ne 0,$ we have already showed the following pointwise estimate (see (\ref{AA1}) and (\ref{EQQ011})).
\begin{equation}\label{eq5}
\big|T(f^{\alpha_1}_1,\ldots,f^{\alpha_m}_m)(x)\big|\lesssim\sum_{j=1}^\infty\prod_{i=1}^m\frac{1}{|2^{j+1}B|}\int_{2^{j+1}B}\big|f_i(y_i)\big|\,dy_i.
\end{equation}
Without loss of generality, we may assume that $p_1=\cdots=p_{\ell}=\min\{p_1,\ldots,p_m\}=1$, and $p_{\ell+1},\ldots,p_m>1$. Using H\"older's inequality, the multiple $A_{\vec{P}}$ condition and Lemma \ref{lem5}, we have the consequences similar to (\ref{EQQ012}) as follows
\begin{equation*}
\begin{split}
&\big|T(f^{\alpha_1}_1,\ldots,f^{\alpha_m}_m)(x)\big|\\
\lesssim&\sum_{j=1}^\infty(\prod_{i=1}^{\ell}\frac{1}{|2^{j+1}B|}\int_{2^{j+1}B}\big|f_i(y_i)\big|\,dy_i)\times\prod_{i=\ell+1}^{m}\frac{1}{|2^{j+1}B|}\int_{2^{j+1}B}\big|f_i(y_i)\big|\,dy_i\\
\lesssim&\sum\limits_{j = 1}^\infty  {\prod\limits_{i = 1}^\ell  {(\frac{1}{{|{2^{j + 1}}B|}}} \int_{{2^{j + 1}}B} | {f_i}({y_i})|{\omega _i}({y_i}){\mkern 1mu} d{y_i}} ){\left( {\mathop {\inf }\limits_{{y_i} \in {2^{j + 1}}B} {\omega _i}({y_i})} \right)^{ - 1}}\\
\times&\prod_{i=\ell+1}^{m}\frac{1}{|2^{j+1}B|}
\left(\int_{2^{j+1}B}\big|f_i(y_i)\big|^{p_i}\omega_i(y_i)\,dy_i\right)^{1/{p_i}}
\left(\int_{2^{j+1}B}\omega_i(y_i)^{1-p'_i}\,dy_i\right)^{1/{p'_i}}\\
\lesssim&(\prod_{i=1}^m \big\|f_i\big\|_{L^{p_i,\kappa}(\omega_i)})
\sum_{j=1}^\infty\nu_{\vec{\omega}}\big(2^{j+1}B\big)^{{(\kappa-1)}/p}.
\end{split}
\end{equation*}
Observe that $\nu_{\vec{\omega}}\in A_{mp}$ with $1\le mp<\infty$. Thus, it follows from the inequality (\ref{AAA5}) that for any $x\in B$,
\begin{align}	
\big|T(f^{\alpha_1}_1,\ldots,f^{\alpha_m}_m)(x)\big|
&\lesssim(\prod_{i=1}^m \big\|f_i\big\|_{L^{p_i,\kappa}(\omega_i)})\cdot\frac{1}{\nu_{\vec{\omega}}(B)^{{(1-\kappa)}/p}}
\sum_{j=1}^\infty\frac{\nu_{\vec{\omega}}(B)^{{(1-\kappa)}/p}}{\nu_{\vec{\omega}}(2^{j+1}B)^{{(1-\kappa)}/p}}\notag\\
&\lesssim(\prod_{i=1}^m \big\|f_i\big\|_{L^{p_i,\kappa}(\omega_i)})\cdot\frac{1}{\nu_{\vec{\omega}}(B)^{{(1-\kappa)}/p}}
\sum_{j=1}^\infty\left(\frac{|B|}{|2^{j+1}B|}\right)^{{\delta(1-\kappa)}/p}\notag\\
&\approx(\prod_{i=1}^m \big\|f_i\big\|_{L^{p_i,\kappa}(\omega_i)})\cdot\frac{1}{\nu_{\vec{\omega}}(B)^{{(1-\kappa)}/p}}.\label{EQQ014}
\end{align}
By using (\ref{EQQ014}), we have 
$${J^{{\alpha _1} \cdots ,{\alpha _m}}} \lesssim \prod\limits_{i = 1}^m {{{\left\| {{f_i}} \right\|}_{{L^{{p_i},\kappa }}({\omega _i})}}} .$$
Combining with (\ref{N4}), we have already finished this proof.
\end{proof}

\subsection{Proofs of Theorem \ref{the2}}
The proofs for the iterated commutator ${g_{\prod {\vec b} }}$ are similar to the commutator ${\mu_{\vec b}}$'s, so we only give the proofs for commutator ${\mu_{\vec b}}$.
\subsubsection{Proof of (i) of Theorem \ref{the2}}
\begin{proof}[Proof:]
For any ball $B=B(x_0,r)$, let $f_i=f^0_i+f^{\infty}_i$, where $f^0_i=f_i\chi_{2B}$, $i=1,\ldots,m$ and $\chi_{2B}$ denotes the characteristic function of $2B$. Then we write
\begin{equation*}
\begin{split}
\prod_{i=1}^m f_i(y_i)&=\prod_{i=1}^m\Big(f^0_i(y_i)+f^{\infty}_i(y_i)\Big)\\
&=\sum_{\alpha_1,\ldots,\alpha_m\in\{0,\infty\}}f^{\alpha_1}_1(y_1)\cdots f^{\alpha_m}_m(y_m)\\
&=\prod\limits_{i = 1}^m {f_i^0} ({y_i}) + \sum\limits_{{\alpha _1} +  \cdots  + {\alpha _m} \ne 0} {f_1^{{\alpha _1}}({y_1}) \cdots f_m^{{\alpha _m}}({y_m})}.
\end{split}
\end{equation*}
Since $T_{\vec b}^N$ is an $m$-sublinear operator, then we have
\begin{equation*}
\begin{split}
&{\nu _{\vec \omega }}{(B)^{ - \frac{\kappa }{p}}}{\left\| {T_{\vec b}^N({f_1}, \cdots ,{f_m})} \right\|_{{L^p}(B,{\nu _{\vec \omega }}dx)}}\\
\le &{\nu _{\vec \omega }}{(B)^{ - \frac{\kappa }{p}}}{\left\| {T_{\vec b}^N(f_1^0, \cdots ,f_m^0)} \right\|_{{L^p}(B,{\nu _{\vec \omega }}dx)}} + \sum\limits_{{\alpha _1} +  \cdots  + {\alpha _m} \ne 0}^{} {{\nu _{\vec \omega }}{{(B)}^{ - \frac{\kappa }{p}}}{{\left\| {T_{\vec b}^N(f_1^{{\alpha _1}}, \cdots ,f_m^{{\alpha _m}})} \right\|}_{_{{L^p}(B,{\nu _{\vec \omega }}dx)}}}} \\
: = &{I^0} + \sum\limits_{{\alpha _1} +  \cdots  + {\alpha _m} \ne 0}^{} {{I^{{\alpha _1}, \cdots ,{\alpha _m}}}}.
\end{split}
\end{equation*}
Due to ${\left\| {{T_{\vec b}}(\vec f)} \right\|_{{L^p}(B,{v_{\vec \omega }}dx)}} \le \sum\limits_{N = 1}^m {{{\left\| {T_{\vec b}^N(\vec f)} \right\|}_{{L^p}(B,{v_{\vec \omega }}dx)}}}$, we merely need to prove:
\begin{equation}
{I^{{\alpha _1}, \cdots ,{\alpha _m}}} \lesssim \prod\limits_{i = 1}^m {{{\left\| {{f_i}} \right\|}_{{L^{{p_i},\kappa }}({\omega _i})}}}, 
\end{equation}
where ${\alpha _i} \in \{ 0,\infty \} ,i = 1, \cdots ,m.$\\
In view of Lemma \ref{lem5}, we have $\nu_{\vec{\omega}}\in A_{mp}$. Applying the boundedness, Lemma \ref{lem1} and Lemma \ref{lem5}, we get
\begin{equation*}
\begin{split}
I^0&\lesssim
\frac{1}{\nu_{\vec{\omega}}(B)^{\kappa/p}}\prod_{i=1}^m\left(\int_{2B}|f_i(x)|^{p_i}\omega_i(x)\,dx\right)^{1/{p_i}}\\
&\lesssim(\prod_{i=1}^m\big\|f_i\big\|_{L^{p_i,\kappa}(\omega_i)})\cdot
\frac{\prod_{i=1}^m \omega_i(2B)^{\kappa/{p_i}}}{\nu_{\vec{\omega}}(B)^{\kappa/p}}\\
&\lesssim(\prod_{i=1}^m\big\|f_i\big\|_{L^{p_i,\kappa}(\omega_i)})\cdot
\frac{\nu_{\vec{\omega}}(2B)^{\kappa/p}}{\nu_{\vec{\omega}}(B)^{\kappa/p}}\\
&\lesssim\prod_{i=1}^m \big\|f_i\big\|_{L^{p_i,\kappa}(\omega_i)}.
\end{split}
\end{equation*}

For $T_{\vec b}^N$, we have
\begin{equation*}
\begin{split}
T_{\vec b}^N(\vec f)(x) &= {\left\| {\int_{\nm} {(Q(x,\vec y))(({b_N}(x) - {b_N}({y_N}))\prod\limits_{j = 1}^m {{f_j}({y_j})} )} d{y_1} \cdots d{y_m}} \right\|_X}\\
&\le \int_{\nm} {{{\left\| {(Q(x,\vec y))(({b_N}(x) - {b_N}({y_N}))\prod\limits_{j = 1}^m {{f_j}({y_j})} )} \right\|}_X}d{y_1} \cdots d{y_m}}\\
&\lesssim \int_{\nm} {\frac{{\left| {(({b_N}(x) - {b_N}({y_N}))\prod\limits_{j = 1}^m {{f_j}({y_j})} } \right|}}{{{{(\sum\limits_{i = 1}^m {\left| {x - {y_i}} \right|} )}^{mn}}}}d{y_1} \cdots d{y_m}} 
\end{split}
\end{equation*}

For the other terms, we first consider the case when $\alpha_1=\cdots=\alpha_m=\infty$. For $x\in B$, we have
\begin{align*}
&\big|T_{\vec b}^N(f^\infty_1,\ldots,f^\infty_m)(x)\big|\\
\lesssim &\int_{{{({\rn})}^m}\backslash {{(2B)}^m}} {\frac{{\left| {{b_N}(x) - {b_N}({y_N})} \right| \cdot |{f_1}({y_1}) \cdots {f_m}({y_m})|}}{{{{(|x - {y_1}| +  \cdots  + |x - {y_m}|)}^{mn}}}}{d\vec y}}\\
\le &\left| {{b_N}(x) - {({b_N})_B}} \right| \cdot \int_{{{({\rn})}^m}\backslash {{(2B)}^m}} {\frac{{|{f_1}({y_1}) \cdots {f_m}({y_m})|}}{{{{(|x - {y_1}| +  \cdots  + |x - {y_m}|)}^{mn}}}}{d\vec y}}\\
+ &\int_{{{({\rn})}^m}\backslash {{(2B)}^m}} {\frac{{\left| {{b_N}({y_N}) - {({b_N})_B}} \right| \cdot |{f_1}({y_1}) \cdots {f_m}({y_m})|}}{{{{(|x - {y_1}| +  \cdots  + |x - {y_m}|)}^{mn}}}}{d\vec y}}\\
\lesssim &\left| {{b_N}(x) - {({b_N})_B}} \right| \cdot \sum\limits_{j = 1}^{\infty} {\frac{1}{{{{\left| {{2^{j + 1}}B} \right|}^m}}}\prod\limits_{i = 1}^m {\int_{{2^{j + 1}}B} {\left| {{f_i}({y_i})} \right|} d{y_i}} }\\
+ &\sum\limits_{j = 1}^{\infty} {\frac{1}{{{{\left| {{2^{j + 1}}B} \right|}^m}}}{\int _{{2^{j + 1}}B}}\left| {{b_N}({y_N}) - {({b_N})_{{2^{j + 1}}B}}} \right| \cdot \left| {{f_N}({y_N})} \right|d{y_N} \cdot \prod\limits_{i \ne N}^{} {{\int _{{2^{j + 1}}B}}\left| {{f_i}({y_i})} \right|d{y_i}} }\\
+ &\sum\limits_{j = 1}^{\infty} {\frac{1}{{{{\left| {{2^{j + 1}}B} \right|}^m}}}{\int _{{2^{j + 1}}B}}\left| {{({b_N})_{{2^{j + 1}}B}} - {({b_N})_B}} \right| \cdot \left| {{f_N}({y_N})} \right|d{y_N} \cdot \prod\limits_{i \ne N}^{} {{\int _{{2^{j + 1}}B}}\left| {{f_i}({y_i})} \right|d{y_i}} }\\
: = &{M_1} + {M_2}+ {M_3}.
\end{align*}
For ${M_1}$:
\begin{align*}
{M_1} &\lesssim \left| {b_N(x) - {({b_N})_B}} \right|\sum\limits_{j = 1}^{\infty} {\frac{1}{{{{\left| {{2^{j + 1}}B} \right|}^m}}}{{\prod\limits_{i = 1}^m {(\int_{{2^{j + 1}}B} {{{\left| {{f_i}({y_i})} \right|}^{{p_i}}}{\omega _i}({y_i})d{y_i})} } }^{\frac{1}{{{p_i}}}}} \cdot (\int_{{2^{j + 1}}B} {{\omega _i}{{({y_i})}^{ - \frac{{{p_i}^\prime }}{{{p_i}}}}}d{y_i}{)^{\frac{1}{{{p_i}^\prime }}}}} }\\
&\lesssim \left| {b_N(x) - {({b_N})_B}} \right|\sum\limits_{j = 1}^\infty  {{v_{\vec \omega }}{{({2^{j + 1}}B)}^{\frac{{\kappa - 1}}{p}}}\prod\limits_{i = 1}^m {{{\left\| {{f_i}} \right\|}_{{L^{{p_i},\kappa }}({\omega _i})}}} }\\
&\lesssim \left| {b_N(x) - {({b_N})_B}} \right|{v_{\vec \omega }}{(B)^{\frac{{\kappa - 1}}{p}}}\prod\limits_{i = 1}^m {{{\left\| {{f_i}} \right\|}_{{L^{{p_i},\kappa }}({\omega _i})}}}.
\end{align*}
where we use this fact:
\begin{equation}
\sum\limits_{j = 1}^\infty  {{v_{\vec \omega }}{{({2^{j + 1}}B)}^{\frac{{\kappa - 1}}{p}}}}  = {v_{\vec \omega }}{(B)^{\frac{{\kappa - 1}}{p}}}\sum\limits_{j = 1}^\infty  {\frac{{{v_{\vec \omega }}{{(B)}^{\frac{{1 - \kappa}}{p}}}}}{{{v_{\vec \omega }}{{({2^{j + 1}}B)}^{\frac{{1 - \kappa}}{p}}}}}}  \lesssim {v_{\vec \omega }}{(B)^{\frac{{\kappa - 1}}{p}}}.
\end{equation}
Then, we have
\begin{align*}
{v_{\vec \omega }}{(B)^{ - \frac{\kappa}{p}}}{\left\| {{M_1}} \right\|_{{L^p}(B,{v_{\vec \omega }}dx)}}&\lesssim (\prod\limits_{i = 1}^m {{{\left\| {{f_i}} \right\|}_{{L^{p_i,\kappa }}({\omega _i})}}}) {(\frac{1}{{{v_{\vec \omega }}(B)}}{\int_B^{} {\left| {b_N(x) - {(b_N)_B}} \right|} ^p}{v_{\vec \omega }}(x)dx)^{\frac{1}{p}}}\\
&\le (\prod\limits_{i = 1}^m {{{\left\| {{f_i}} \right\|}_{{L^{p_i,\kappa }}({\omega _i})}}} ){v_{\vec \omega }}{(B)^{ - \frac{1}{p}}}{(\int_B^{} {{{\left| {{b_N}(x) - {{({b_N})}_B}} \right|}^{pr'}}} dx)^{\frac{1}{{pr'}}}}{(\int_B^{} {{v_{\vec \omega }}^r} )^{\frac{1}{{pr}}}}\\
&\lesssim {\left\| b_N \right\|_{BMO}}\prod\limits_{i = 1}^m {{{\left\| {{f_i}} \right\|}_{{L^{p_i,\kappa }}({\omega _i})}}},
\end{align*}
where the last inequality is valid, since we use Lemma \ref{lem2} and reverse H\"older's inequality.
For ${M_3}$:
\begin{align*}
{M_3} &\lesssim \sum\limits_{j = 1}^{\infty} {\frac{{j \cdot {{\left\| b_N \right\|}_{BMO}}}}{{{{\left| {{2^{j + 1}}B} \right|}^m}}}\prod\limits_{i = 1}^m {\int_{{2^{j + 1}}B} {\left| {{f_i}({y_i})} \right|} d{y_i}} }\\
&\lesssim {\left\| b_N \right\|_{BMO}} {v_{\vec \omega }}{(B)^{\frac{{\kappa - 1}}{p}}} \prod\limits_{j = 1}^m {{{\left\| {{f_i}} \right\|}_{{L^{p_i,\kappa }}({\omega _i})}}},
\end{align*}
where we do a simple calculation, and use the following fact, see (\cite{Torchinsky}, p.206),
\begin{equation}\label{eq7}
\left| {{(b_N)_{{2^{j + 1}}B}} - {(b_N)_B}} \right| \lesssim ({j+1}){\left\| b_N \right\|_{BMO}}.
\end{equation}
Then, we have
$${v_{\vec \omega }}{(B)^{ - \frac{\kappa}{p}}}{\left\| {{M_3}} \right\|_{{L^p}(B,{v_{\vec \omega }}dx)}} \lesssim {\left\| b_N \right\|_{BMO}}\prod\limits_{i = 1}^m {{{\left\| {{f_i}} \right\|}_{{L^{p_i,\kappa }}({\omega _i})}}}.$$
For ${M_2}$:
\begin{equation}\label{EQQ018}
\begin{aligned}
&\int_{{2^{j + 1}}B}^{} {\left| {(b_N)(y_N) - {(b_N)_{{2^{j + 1}}B}}} \right|\left| {f(y_N)} \right|dy_N}\\
&\le {\left[ {\int_{{2^{j + 1}}B}^{} {{{\left| {{b_N}({y_N}) - {{({b_N})}_{{2^{j + 1}}B}}} \right|}^{{p_N}^\prime }}{\omega _N}{{({y_N})}^{ - \frac{{{p_N}^\prime }}{{{p_N}}}}}d{y_N}} } \right]^{\frac{1}{{{p_N}^\prime }}}} \cdot {(\int_{{2^{j + 1}}B}^{} {{{\left| {{f_N}} \right|}^{{p_N}}}} {\omega _N})^{\frac{1}{{{p_N}}}}}\\
&\le {\left\| {{f_N}} \right\|_{{L^{{p_N},\kappa}}({\omega _N})}}{\left[ {\int_{{2^{j + 1}}B}^{} {{{\left| {b_N({y_N}) - {(b_N)_{{2^{j + 1}}B}}} \right|}^{{{p_N}^\prime}}}{\omega _N}{{({y_N})}^{ - \frac{{{{p_N}^\prime}}}{{{p_N}}}}}d{y_N}} } \right]^{\frac{1}{{{{p_N}^\prime}}}}} \cdot {\omega _N}{({2^{j + 1}}B)^{\frac{\kappa}{{{p_N}}}}}.
\end{aligned}
\end{equation}
Now, let us prove that the following result is valid:
\begin{equation}\label{EQQ019}
{\left[ {\int_{{2^{j + 1}}B}^{} {{{\left| {b_N({y_N}) - {(b_N)_{{2^{j + 1}}B}}} \right|}^{{{p_N}^\prime}}}{\omega _N}{{({y_N})}^{ - \frac{{{{p_N}^\prime}}}{{{p_N}}}}}d{y_N}} } \right]^{\frac{1}{{{{p_N}^\prime}}}}} \lesssim {\left\| b_N \right\|_{BMO}}{v_N}{({2^{j + 1}}B)^{\frac{1}{{{{p_N}^\prime}}}}}.
\end{equation}
In fact, there exists a $\theta  > 1$, such that ${\omega _N}^{ - \frac{{{{p_N}^\prime}}}{{{p_N}}}}: = {v_N} \in R{H_\theta }$, since ${v_N} \in {A_{{{p_N}^\prime}}}$.
Then, we have:$${(\frac{1}{{\left| {{2^{j + 1}}B} \right|}}\int_{{2^{j + 1}}B}^{} {{{v_N}^\theta }} )^{\frac{1}{\theta }}} \lesssim (\frac{1}{{\left| {{2^{j + 1}}B} \right|}}\int_{{2^{j + 1}}B}^{} {v_N} ).$$ 
Thus, we have
\begin{equation*}
\begin{aligned}
&{\left[ {\int_{{2^{j + 1}}B}^{} {{{\left| {b_N({y_N}) - {(b_N)_{{2^{j + 1}}B}}} \right|}^{{{p_N}^\prime}}}{\omega _N}{{({y_N})}^{ - \frac{{{{p_N}^\prime}}}{{{p_N}}}}}d{y_N}} } \right]^{\frac{1}{{{{p_N}^\prime}}}}}\\
&\lesssim {\left[ {\int_{{2^{j + 1}}B}^{} {{{\left| {b_N({y_N}) - {(b_N)_{{2^{j + 1}}B}}} \right|}^{{{p_N}^\prime}\theta '}}d{y_N}} } \right]^{\frac{1}{{{{p_N}^\prime}\theta '}}}} \cdot {\left| {{2^{j + 1}}B} \right|^{\frac{1}{{{{p_N}^\prime}\theta '}}}} \cdot {v_N}{({2^{j + 1}}B)^{\frac{1}{{{{p_N}^\prime}}}}}\\
&\lesssim {\left\| b_N \right\|_{BMO}}{v_N}{({2^{j + 1}}B)^{\frac{1}{{{{p_N}^\prime}}}}}.
\end{aligned}
\end{equation*}
We have proved (\ref{EQQ019}) is valid, and then, we apply (\ref{EQQ019}) back to the proof in (\ref{EQQ018}):
\begin{equation*}
\begin{aligned}
&\int_{{2^{j + 1}}B}^{} {\left| {{b_N}({y_N}) - {{({b_N})}_{{2^{j + 1}}B}}} \right|\left| {{f_N}({y_N})} \right|d{y_N}}\\
\lesssim &{\left\| {{b_N}} \right\|_{BMO}}{\left\| {{f_N}} \right\|_{{L^{{p_N},\kappa }}({\omega _N})}}{v_N}{({2^{j + 1}}B)^{\frac{1}{{{p_N}^\prime }}}} \cdot {\omega _N}{({2^{j + 1}}B)^{\frac{\kappa }{{{p_N}}}}}.
\end{aligned}
\end{equation*}
We get the boundedness of ${M_2}$:
\begin{align*}
{M_2}&\lesssim \sum\limits_{j = 1}^\infty  {\frac{1}{{{{\left| {{2^{j + 1}}B} \right|}^{m - 1}}}}(} \prod\limits_{i \ne N} {\int_{{2^{j + 1}}B} {\left| {{f_i}({y_i})} \right|d{y_i}} } ) \cdot {\left\| {{b_N}} \right\|_{BMO}}{\left\| {{f_N}} \right\|_{{L^{{p_N},\kappa }}({\omega _N})}}{v_N}{({2^{j + 1}}B)^{\frac{1}{{{p_N}^\prime }}}} \cdot {\omega _N}{({2^{j + 1}}B)^{\frac{\kappa }{{{p_N}}}}}\\
&\lesssim {\left\| {{b_N}} \right\|_{BMO}}(\prod\limits_{i = 1}^m {{{\left\| {{f_i}} \right\|}_{{L^{{p_i},\kappa }}({\omega _N})}}} )\sum\limits_{j = 1}^\infty  {\frac{1}{{{{\left| {{2^{j + 1}}B} \right|}^m}}} \cdot \prod\limits_{i = 1}^m {{\omega _i}{{({2^{j + 1}}B)}^{\frac{{\kappa  - 1}}{{{p_i}}}}}{v_i}{{({2^{j + 1}}B)}^{\frac{1}{{{p_i}^\prime }}}}} }\\
&\lesssim {\left\| b_N \right\|_{BMO}}(\prod\limits_{i = 1}^m {{{\left\| {{f_i}} \right\|}_{{L^{{p_i},\kappa }}({\omega _i})}}}) \sum\limits_{j = 1}^\infty  {{v_{\vec \omega }}{{({2^{j + 1}}B)}^{\frac{{\kappa - 1}}{p}}}}\\
&\lesssim {\left\| b_N \right\|_{BMO}} {v_{\vec \omega }}{(B)^{\frac{{\kappa - 1}}{p}}} \prod\limits_{i = 1}^m {{{\left\| {{f_i}} \right\|}_{{L^{p_i,\kappa }}({\omega _i})}}},
\end{align*}
where we set ${v_i}={\omega _i}^{ - \frac{{{{p_i}^\prime}}}{{{p_i}}}},$ some of the details are similar to the previous proof which we omit here, and then we have 
$${v_{\vec \omega }}{(B)^{ - \frac{\kappa}{p}}}{\left\| {{M_2}} \right\|_{{L^p}(B,{v_{\vec \omega }}dx)}} \lesssim {\left\| b_N \right\|_{BMO}}\prod\limits_{i = 1}^m {{{\left\| {{f_i}} \right\|}_{{L^{{p_i},\kappa }}({\omega _i})}}}.$$
Consequently, we have
\begin{equation}
{I^{\infty , \cdots ,\infty }} \lesssim {\left\| b_N \right\|_{BMO}}\prod\limits_{i = 1}^m {{{\left\| {{f_i}} \right\|}_{{L^{p_i,\kappa }}({\omega _i})}}}.
\end{equation}
Without loss of generality, we may assume that ${\alpha _1} =  \cdots  = {\alpha _\ell } = \infty$, and ${\alpha _{l + 1}} =  \cdots  = {\alpha _m} = 0.$
If $N \in \{ 1, \cdots ,l\},$ we have 
\begin{align*}
&|T_{\vec b}^N(f_1^\infty , \ldots ,f_\ell ^\infty ,f_{\ell  + 1}^0, \ldots ,f_m^0)(x)|\\
\lesssim&\int_{{{[{{(2B)}^c}]}^l}}^{} {\int_{{{(2B)}^{m - l}}}^{} {\frac{{\left| {b_N(x) - (b_N)({y_N})} \right| \cdot |{f_1}({y_1}) \cdots {f_m}({y_m})|}}{{{{(|x - {y_1}| +  \cdots  + |x - {y_m}|)}^{mn}}}}d\vec y} }\\
\lesssim&\left| {b_N(x) - {(b_N)_B}} \right|\sum\limits_{j = 1}^{\infty} {\frac{1}{{{{\left| {{2^{j + 1}}B} \right|}^m}}}(\prod\limits_{i = 1}^m {\int_{{2^{j + 1}}B}^{} {\left| {{f_i}({y_i})} \right|} } d{y_i})}\\
+& \sum\limits_{j = 1}^{\infty} {\frac{1}{{{{\left| {{2^{j + 1}}B} \right|}^m}}}(\prod\limits_{i \ne N}^{} {\int_{{2^{j + 1}}B}^{} {\left| {{f_i}({y_i})} \right|} } d{y_i}) \cdot \int_{{2^{j + 1}}B} {\left| {b_N({y_N}) - {(b_N)_{{2^{j + 1}}B}}} \right|} \left| {{f_N}({y_N})} \right|} d{y_N}\\
+ &\sum\limits_{j = 1}^{\infty} {\frac{1}{{{{\left| {{2^{j + 1}}B} \right|}^m}}}(\prod\limits_{i \ne N}^{} {\int_{{2^{j + 1}}B}^{} {\left| {{f_i}({y_i})} \right|} } d{y_i}) \cdot \int_{{2^{j + 1}}B} {\left| {{(b_N)_{{2^{j + 1}}B}} - {(b_N)_B}} \right|} \left| {{f_N}({y_N})} \right|} d{y_N}\\
= & {M_1}+{M_2}+{M_3}.
\end{align*}
We just consider the last case for now. For any $x\in B$, if $N \in \{ l + 1, \cdots ,m\},$ we have 
\begin{align*}
&|T_{\vec b}^N(f_1^\infty , \ldots ,f_\ell ^\infty ,f_{\ell  + 1}^0, \ldots ,f_m^0)(x)|\\
\lesssim &\sum\limits_{j = 1}^{\infty} {\int_{{{({2^{j + 1}}B\backslash {2^j}B)}^l}}^{} {\int_{{{(2B)}^{m - l}}}^{} {\frac{{\left| {b_N(x) - b_N({y_N})} \right| \cdot |{f_1}({y_1}) \cdots {f_m}({y_m})|}}{{{{(|x - {y_1}| +  \cdots  + |x - {y_m}|)}^{mn}}}}d\vec y} } }\\
\lesssim &\left| {b_N(x) - {(b_N)_B}} \right|\sum\limits_{j = 1}^{\infty} {\frac{1}{{{{\left| {{2^{j + 1}}B} \right|}^m}}}\int_{{2^{j + 1}}B} {|{f_N}({y_N})|} d{y_N}\prod\limits_{i \ne N}^{} {\int_{{2^{j + 1}}B} {\left| {{f_i}({y_i})} \right|} d{y_i}} }\\
+ &\sum\limits_{j = 1}^{\infty} {\frac{1}{{{{\left| {{2^{j + 1}}B} \right|}^m}}}\int_{{2^{j + 1}}B} {\left| {b_N({y_N}) - {(b_N)_{{2^{j + 1}}B}}} \right| \cdot |{f_N}({y_N})|} d{y_N}\prod\limits_{i \ne N}^{} {\int_{{2^{j + 1}}B} {\left| {{f_i}({y_i})} \right|} d{y_i}} }\\
+ &\sum\limits_{j = 1}^{\infty} {\frac{1}{{{{\left| {{2^{j + 1}}B} \right|}^m}}}\int_{{2^{j + 1}}B} {\left| {{(b_N)_{{2^{j + 1}}B}} - {(b_N)_B}} \right| \cdot |{f_N}({y_N})|} d{y_N}\prod\limits_{i \ne N}^{} {\int_{{2^{j + 1}}B} {\left| {{f_i}({y_i})} \right|} d{y_i}} }\\
= & {M_1}+{M_2}+{M_3}.
\end{align*}
Since ${M_i}$ has the boundedness that we need, for $i=1,2,3$, combining all of the above results, we finish this proof of (i).
\end{proof}

\subsubsection{Proof of (ii) of Theorem \ref{the2}}
\begin{proof}[Proof:]
For any ball $B=B(x_0,r)$, decompose $f_i=f^0_i+f^{\infty}_i$, where $f^0_i=f_i\chi_{2B}$, $i=1,\ldots,m$. For each $\lambda>0$, we have
\begin{equation*}
\begin{split}
&{\nu _{\vec \omega }}{(B)^{ - \frac{\kappa }{p}}}{\left\| {T_{\vec b}^N({f_1}, \cdots ,{f_m})} \right\|_{W{L^p}(B,{\nu _{\vec \omega }}dx)}}\\
\lesssim&{\nu _{\vec \omega }}{(B)^{ - \frac{\kappa }{p}}}{\left\| {T_{\vec b}^N(f_1^0, \cdots ,f_m^0)} \right\|_{W{L^p}(B,{\nu _{\vec \omega }}dx)}} + \sum\limits_{{\alpha _1} +  \cdots  + {\alpha _m} \ne 0}^{} {{\nu _{\vec \omega }}{{(B)}^{ - \frac{\kappa }{p}}}{{\left\| {T_{\vec b}^N(f_1^{{\alpha _1}}, \cdots ,f_m^{{\alpha _m}})} \right\|}_{_{W{L^p}(B,{\nu _{\vec \omega }}dx)}}}}\\
: = &{J^0} + \sum\limits_{{\alpha _1} +  \cdots  + {\alpha _m} \ne 0} {{J^{{\alpha _1}, \cdots ,{\alpha _m}}}}.
\end{split}
\end{equation*}
Due to ${\left\| {{T_{\vec b}}(\vec f)} \right\|_{W{L^p}(B,{v_{\vec \omega }}dx)}}\lesssim \sum\limits_{N = 1}^m {{{\left\| {T_{\vec b}^N(\vec f)} \right\|}_{W{L^p}(B,{v_{\vec \omega }}dx)}}}$, we merely need to prove:
\begin{equation}\label{EQQ013}
{J^{{\alpha _1}, \cdots ,{\alpha _m}}} \lesssim \prod_{i=1}^m \big\|f_i\big\|_{L^{p_i,\kappa}(\omega_i)}.
\end{equation}
where ${\alpha _i} \in \{ 0,\infty \} ,i = 1, \cdots ,m.$\\
In view of Lemma \ref{lem4}, we have $\nu_{\vec{\omega}}\in A_{mp}$. Applying the boundedness, Lemma \ref{lem1} and Lemma \ref{lem5}, we get
\begin{equation*}
J^0\lesssim{\nu _{\vec \omega }}{(B)^{ - \frac{\kappa }{p}}}{\prod\limits_{i = 1}^m {\left\| {{f_i}} \right\|} _{{L^{{p_i}}}(2B,{\omega_i}dx)}}\\
\lesssim\frac{{\prod\limits_{i = 1}^m {{\omega_i}} {{(2B)}^{\kappa /{p_i}}}}}{{{\nu _{\vec \omega }}{{(B)}^{\kappa /p}}}}{\prod\limits_{i = 1}^m {\left\| {{f_i}} \right\|} _{{L^{{p_i},\kappa }}({\omega _i})}}\lesssim\prod_{i=1}^m \big\|f_i\big\|_{L^{p_i,\kappa}(\omega_i)}.
\end{equation*}
Thus, we have ${J^0} \lesssim \prod_{i=1}^m \big\|f_i\big\|_{L^{p_i,\kappa}(\omega_i)}.$\\
In the proof of (i), we have already showed the following important estimates. For ${\alpha _1} +  \cdots  + {\alpha _m} \ne 0,$ we have
\begin{equation}
|T_{\vec b}^N(f_1^{{\alpha _1}}, \ldots ,f_m^{{\alpha _m}})(x)|\lesssim {M_1} + {M_2} + {M_3};
\end{equation}
\begin{equation}\label{AAA2}
{M_1} \lesssim \left| {b_N(x) - {(b_N)_B}} \right|{v_{\vec \omega }}{(B)^{\frac{{\kappa - 1}}{p}}}\prod\limits_{i = 1}^m {{{\left\| {{f_i}} \right\|}_{{L^{{p_i},\kappa }}({\omega _i})}}};
\end{equation}
\begin{equation}\label{AAA3}
{M_2} \lesssim {\left\| b_N \right\|_{BMO}} {v_{\vec \omega }}{(B)^{\frac{{\kappa - 1}}{p}}} \prod\limits_{i = 1}^m {{{\left\| {{f_i}} \right\|}_{{L^{p_i,\kappa }}({\omega _i})}}};
\end{equation}
\begin{equation}\label{AAA4}
{M_3} \lesssim {\left\| b_N \right\|_{BMO}} {v_{\vec \omega }}{(B)^{\frac{{\kappa - 1}}{p}}} \prod\limits_{i = 1}^m {{{\left\| {{f_i}} \right\|}_{{L^{p_i,\kappa }}({\omega _i})}}},
\end{equation}
where for some $p_i =1,$ we can still get (\ref{AAA2}), (\ref{AAA3}) and (\ref{AAA4}), since their proof is similar to before.\\
Obviously, we have
$${J^{{\alpha _1}, \cdots ,{\alpha _m}}} \lesssim \sum\limits_{i = 1}^3 {{v_{\vec \omega }}{{(B)}^{ - \frac{\kappa }{p}}}{{\left\| {{M_i}} \right\|}_{{L^p}(B,{v_{\vec \omega }}dx)}}}  \lesssim {\left\| {{b_N}} \right\|_{BMO}}\prod\limits_{i = 1}^m {{{\left\| {{f_i}} \right\|}_{{L^{{p_i},\kappa }}({\omega _i})}}}.$$ 
Consequently, we finish the proof of (ii).
\end{proof}

\subsection{Proof of Theorem \ref{thm3}}
The proof for the iterated commutator ${T_{\prod {\vec b}}}$ is also similar to the commutator ${T_{\vec b}}$'s, so we also only give the proof for commutator ${T_{\vec b}}$.
\begin{proof}[Proof:]
For any ball $B=B(x_0,r)$, decompose $f_i=f^0_i+f^{\infty}_i$, where $f^0_i=f_i\chi_{2B}$, $i=1,\ldots,m$. Here, we merely think about the follows commutator, due to the properties of ${T_{\vec b}}$.
	\begin{equation*}
T_b^1(\vec f)(x) = {\left\| {\int_{{\rm{\nm}}} {(Q(x,\vec y))((b(x) - b({y_1}))\prod\limits_{j = 1}^m {{f_j}({y_j})} )d{y_1} \cdots d{y_m}} } \right\|_X}.
	\end{equation*}
	According to inequalities (see \cite{Gra1}, p.12, exercises 1.1.4.), we have
	\begin{align*}
	&{\nu_{\vec{\omega}}(B)^{-m\kappa}}{[{\nu _{\vec \omega}}(\{ x \in B:|T_{b}^1(\vec f)(x)| > {\lambda ^m}\} )]^m}\notag\\
	\lesssim& {\nu_{\vec{\omega}}(B)^{-m\kappa}}
{[{\nu _{\vec \omega}}(\{ x \in B:|T_{b}^1(f_1^0, \ldots ,f_m^0)(x)| > {\lambda ^m}/{2^m}\} )]^m}\notag\\
	+&\sum\limits_{({\alpha_1}, \ldots ,{\alpha _m}) \ne 0} {\nu_{\vec{\omega}}(B)^{-m\kappa}}
	{[{\nu _{\vec \omega}}(\{ x \in B:|T_{b}^1(f_1^{{\alpha _1}}, \ldots ,f_m^{{\alpha _m}})(x)| > {\lambda ^m}/{2^m}\} )]^m}\notag\\
	:=&H^{0}+\sum\limits_{({\alpha_1}, \ldots ,{\alpha _m}) \ne 0} H^{\alpha_1,\dots,\alpha_m}.
	\end{align*}
	Note that $\Phi(t)=t(1+\log^+t)$ satisfies the following condition, see also (\cite{Gra2}, p.197): for $C>1$, for every $t>0$, $$\Phi (\frac{t}{C}) \le \frac{{\Phi (t)}}{C}.$$
	Thus, combining (\ref{eq3}) and above inequality, we deduce
	\begin{equation*}
	\begin{split}
	H^{0}&\lesssim {\nu_{\vec{\omega}}(B)^{-m\kappa}}
	\prod_{i=1}^m\bigg(\int_{\mathbb R^n}\Phi\bigg(\frac{2|f^0_i(x)|}{\lambda}\bigg)\cdot \omega_i(x)\,dx\bigg)\\
	&\lesssim {\nu_{\vec{\omega}}(B)^{-m\kappa}}
	\prod_{i=1}^m\bigg(\int_{2B}\Phi\bigg(\frac{|f_i(x)|}{\lambda}\bigg)\cdot \omega_i(x)\,dx\bigg)\\
	&\leq {\nu_{\vec{\omega}}(B)^{-m\kappa}}
	\prod_{i=1}^m \omega_i(2B)\cdot\bigg\|\Phi\bigg(\frac{|f_i|}{\lambda}\bigg)\bigg\|_{L\log L(\omega_i),2B},
	\end{split}
	\end{equation*}
	where last inequality is valid, since we have used the estimate \eqref{main esti1}. Due to  Lemma \ref{lem4}, we can see $\nu_{\vec{\omega}}\in A_1$ and $\omega_i^{1/{m}}\in A_1$ ($i=1,2,\dots,m$). Thus, using \eqref{N1}, we have
	\begin{equation*}
	\begin{split}
	H^{0}
	\lesssim\prod_{i=1}^m\bigg\|\Phi\bigg(\frac{|f_i|}{\lambda}\bigg)\bigg\|_{(L\log L)^{1,\kappa}(\omega_i)}
	\frac{1}{\nu_{\vec{\omega}}(B)^{m\kappa}}
	\prod_{i=1}^m \omega_i(2B)^{\kappa}\lesssim\prod_{i=1}^m\bigg\|\Phi\bigg(\frac{|f_i|}{\lambda}\bigg)\bigg\|_{(L\log L)^{1,\kappa}(\omega_i)}.
	\end{split}
	\end{equation*}
	Now, we consider $H^{\alpha_1,\dots,\alpha_m}$ for $({\alpha _1}, \ldots ,{\alpha _m}) \ne 0$. For any $x\in B$,
	\begin{equation*}
	\begin{split}
T_b^1(\vec f)(x) \le& {\left\| {\int_{{\rm{\nm}}} {(Q(x,\vec y))((b(x) - {b_B})\prod\limits_{j = 1}^m {{f_j}({y_j})} )d{y_1} \cdots d{y_m}} } \right\|_X}\\
+& {\left\| {\int_{{\rm{\nm}}} {(Q(x,\vec y))((b(y_1) - {b_B})\prod\limits_{j = 1}^m {{f_j}({y_j})} )d{y_1} \cdots d{y_m}} } \right\|_X}\\
:=&L(\vec f)(x) + \widetilde L(\vec f)(x).
	\end{split}
	\end{equation*}
	So we have,
	\begin{equation*}
	\begin{split}
	H^{\alpha_1,\dots,\alpha_m}
	\lesssim &{\nu_{\vec{\omega}}(B)^{-m\kappa}}
	\Big[\nu_{\vec{\omega}}\Big(\Big\{x\in B:L(\vec f)(x)>\lambda^m/{2^{m+1}}\Big\}\Big)\Big]^m\\
	+&{\nu_{\vec{\omega}}(B)^{-m\kappa}}
	\Big[\nu_{\vec{\omega}}\Big(\Big\{x\in B:\widetilde L(\vec f)(x)>\lambda^m/{2^{m+1}}\Big\}\Big)\Big]^m\\
	:=&{L^{{\alpha_1}, \ldots ,{\alpha_m}}} + {\widetilde L^{{\alpha_1}, \ldots ,{\alpha _m}}}.
	\end{split}
	\end{equation*}
	Combining (\ref{eq5}), we get
	\begin{equation*}
	\begin{split}
	{L^{{\alpha _1}, \ldots ,{\alpha _m}}}&\lesssim {\nu_{\vec{\omega}}(B)^{-m\kappa}}
	\frac{2^{m+1}}{\lambda^m}
	\bigg(\int_{B}\big|L(f^{\alpha_1}_1,f^{\alpha_2}_2,\ldots,f^{\alpha_m}_m)(x)\big|^{\frac{\,1\,}{m}}\nu_{\vec{\omega}}(x)\,dx\bigg)^m\\
	&\lesssim {\nu_{\vec{\omega}}(B)^{-m\kappa}}
	\sum_{j=1}^\infty\bigg(\prod_{i=1}^m\frac{1}{|2^{j+1}B|}
	\int_{2^{j+1}B}\frac{|f_i(y_i)|}{\lambda}\,dy_i\bigg)\bigg(\int_B\big|b(x)-b_{B}\big|^{\frac{\,1\,}{m}}\nu_{\vec{\omega}}(x)\,dx\bigg)^m.
	\end{split}
	\end{equation*}
Since $\nu_{\vec{w}}\in A_1$, there exists $h \in (1,\infty )$, such that $\nu_{\vec{\omega}}\in RH_h$. Then, we have
\begin{equation*}
\begin{split}
\int_B\big|b(x)-b_{B}\big|^{\frac{\,1\,}{m}}\nu_{\vec{\omega}}(x)\,dx
&\leq|B|\bigg(\frac{1}{|B|}\int_B\big|b(x)-b_{B}\big|^{h'/m}\,dx\bigg)^{1/h'}\bigg(\frac{1}{|B|}\int_B\nu_{\vec{\omega}}(x)^h\,dx\bigg)^{1/h}\\
&\lesssim{\nu _{\vec \omega}}{(B)^{\frac{1}{m}}}\left\| b \right\|_{BMO}^{\frac{1}{m}}.
\end{split}
\end{equation*}
where the last inequality follows from inequality (\ref{eq6}) and reverse H\"{o}lder inequality.
	\begin{equation*}
	{L^{{\alpha _1}, \ldots ,{\alpha _m}}}
	\lesssim\|b\|_{BMO}\cdot\nu_{\vec{\omega}}(B)^{m(1-\kappa)}
	\sum_{j=1}^\infty\bigg(\prod_{i=1}^m\frac{1}{|2^{j+1}B|}
	\int_{2^{j+1}B}\frac{|f_i(y_i)|}{\lambda}\,dy_i\bigg).
	\end{equation*}
	Back to ${L^{{\alpha _1}, \ldots ,{\alpha _m}}}$, we have the estimates as follows.
	\begin{equation*}
	\begin{split}
	&{L^{{\alpha _1}, \ldots ,{\alpha _m}}}\\
	\lesssim&\|b\|_{BMO}\cdot\nu_{\vec{\omega}}(B)^{m(1-\kappa)}\sum_{j=1}^\infty\prod_{i=1}^m\bigg(\frac{1}{|2^{j+1}B|}
	\int_{2^{j+1}B}\frac{|f_i(y_i)|}{\lambda}\cdot \omega_i(y_i)\,dy_i\bigg)
	\left(\inf_{y_i\in 2^{j+1}B}\omega_i(y_i)\right)^{-1}\\
	\lesssim&\|b\|_{BMO}\cdot\nu_{\vec{\omega}}(B)^{m(1-\kappa)}
	\sum_{j=1}^\infty\frac{1}{\nu_{\vec{\omega}}(2^{j+1}B)^m}
	\prod_{i=1}^m\int_{2^{j+1}B}\Phi\bigg(\frac{|f_i(y_i)|}{\lambda}\bigg)\cdot \omega_i(y_i)\,dy_i\\
	\lesssim&\|b\|_{BMO}\cdot\nu_{\vec{\omega}}(B)^{m(1-\kappa)}
	\sum_{j=1}^\infty\frac{1}{\nu_{\vec{\omega}}(2^{j+1}B)^m}
	\prod_{i=1}^m\omega_i\big(2^{j+1}B\big)\bigg\|\Phi\bigg(\frac{|f_i|}{\lambda}\bigg)\bigg\|_{L\log L(\omega_i),2^{j+1}B},
	\end{split}
	\end{equation*}
	where the last inequality follows from the previous estimate \eqref{main esti1}. In view of \eqref{N1}, the last expression is bounded by
	\begin{equation*}
	\begin{split}
	&\|b\|_{BMO}\cdot\nu_{\vec{\omega}}(B)^{m(1-\kappa)}
	\sum_{j=1}^\infty\frac{1}{\nu_{\vec{\omega}}(2^{j+1}B)^m}\prod_{i=1}^m\bigg\|\Phi\bigg(\frac{|f_i|}{\lambda}\bigg)\bigg\|_{(L\log L)^{1,\kappa}(\omega_i)}
	\prod_{i=1}^m \omega_i(2^{j+1}B)^{\kappa}\\
	\lesssim&\|b\|_{BMO}\prod_{i=1}^m\bigg\|\Phi\bigg(\frac{|f_i|}{\lambda}\bigg)\bigg\|_{(L\log L)^{1,\kappa}(\omega_i)}
	\sum_{j=1}^\infty\frac{\nu_{\vec{\omega}}(B)^{m(1-\kappa)}}{\nu_{\vec{\omega}}(2^{j+1}B)^{m(1-\kappa)}}\\
	\lesssim&\|b\|_{BMO}\prod_{i=1}^m\bigg\|\Phi\bigg(\frac{|f_i|}{\lambda}\bigg)\bigg\|_{(L\log L)^{1,\kappa}(\omega_i)}.
	\end{split}
	\end{equation*}
Applying the pointwise estimates \eqref{eq5}, we have
	\begin{equation*}
	\begin{split}
	{\widetilde L^{{\alpha _1}, \ldots ,{\alpha _m}}}&\lesssim{\nu_{\vec{\omega}}(B)^{-m\kappa}}
	\frac{2^{m+1}}{\lambda^m}
	\bigg(\int_{B}\big|\widetilde L(f^{\alpha_1}_1,f^{\alpha_2}_2,\ldots,f^{\alpha_m}_m)(x)\big|^{\frac{\,1\,}{m}}
	\nu_{\vec{\omega}}(x)\,dx\bigg)^m\\
	&\lesssim \nu_{\vec{\omega}}(B)^{m(1-\kappa)}
	\sum_{j=1}^\infty\bigg(\prod_{i=2}^m\frac{1}{|2^{j+1}B|}\int_{2^{j+1}B}\frac{|f_i(y_i)|}{\lambda}\,dy_i\bigg)\\
	&\bigg(\frac{1}{|2^{j+1}B|}\int_{2^{j+1}B}
	\big|b(y_1)-b_{B}\big|\cdot\frac{|f_1(y_1)|}{\lambda}\,dy_1\bigg)\\
	&\lesssim\nu_{\vec{\omega}}(B)^{m(1-\kappa)}
	\sum_{j=1}^\infty\bigg(\prod_{i=2}^m\frac{1}{|2^{j+1}B|}
	\int_{2^{j+1}B}\frac{|f_i(y_i)|}{\lambda}\omega_i(y_i)\,dy_i\bigg)\\
	&\bigg(\frac{1}{|2^{j+1}B|}\int_{2^{j+1}B}
	\big|b(y_1)-b_{B}\big|\cdot\frac{|f_1(y_1)|}{\lambda}\omega_1(y_1)\,dy_1\bigg)\prod_{i=1}^m\left(\inf_{y_i\in 2^{j+1}B}\omega_i(y_i)\right)^{-1}\\
	&\lesssim \nu_{\vec{\omega}}(B)^{m(1-\kappa)}
	\times\sum_{j=1}^\infty\frac{1}{\nu_{\vec{\omega}}(2^{j+1}B)^m}
	\bigg(\prod_{i=2}^m\int_{2^{j+1}B}\frac{|f_i(y_i)|}{\lambda}\omega_i(y_i)\,dy_i\bigg)\\
	&\bigg(\int_{2^{j+1}B}\big|b(y_1)-b_{B}\big|\cdot\frac{|f_1(y_1)|}{\lambda}\omega_1(y_1)\,dy_1\bigg),
	\end{split}
	\end{equation*}
	where in the last inequality we have used the $A_{(1,\dots,1)}$ condition. Next, we have two estimates as follows.
	\begin{equation*}
	\begin{split}
	\int_{2^{j+1}B}\frac{|f_i(y_i)|}{\lambda}\omega_i(y_i)\,dy_i&\leq\int_{2^{j+1}B}\Phi\bigg(\frac{|f_i(y_i)|}{\lambda}\bigg)\cdot \omega_i(y_i)\,dy_i\\
	&\leq \omega_i\big(2^{j+1}B\big)\bigg\|\Phi\bigg(\frac{|f_i|}{\lambda}\bigg)\bigg\|_{L\log L(w_i),2^{j+1}B}.
	\end{split}
	\end{equation*}
	Using the inequality \eqref{Wholder}, we obtain
	\begin{equation*}
	\begin{split}
	&\int_{2^{j+1}B}\big|b(y_1)-b_{B}\big|\cdot\frac{|f_1(y_1)|}{\lambda}\omega_1(y_1)\,dy_1\\
	\leq&\int_{2^{j+1}B}\big|b(y_1)-b_{B}\big|\cdot\Phi\bigg(\frac{|f_1(y_1)|}{\lambda}\bigg)\omega_1(y_1)\,dy_1\\
	\lesssim &\omega_1\big(2^{j+1}B\big)
	\big\|b-b_{B}\big\|_{\exp L(\omega_1),2^{j+1}B}
	\bigg\|\Phi\bigg(\frac{|f_1|}{\lambda}\bigg)\bigg\|_{L\log L(\omega_1),2^{j+1}B}.
	\end{split}
	\end{equation*}
	Combining the inequality \eqref{eq7}, \eqref{BMOwang} and \eqref{Wholder}, we can get
	\begin{equation*}
	\int_{2^{j+1}B}\big|b(y_1)-b_{B}\big|\cdot\frac{|f_1(y_1)|}{\lambda}\omega_1(y_1)\,dy_1\lesssim (j+1)\|b\|_{BMO}\omega_1\big(2^{j+1}B\big)
	\bigg\|\Phi\bigg(\frac{|f_1|}{\lambda}\bigg)\bigg\|_{L\log L(\omega_1),2^{j+1}B}.
	\end{equation*}
	Thus, we can deduce the results from the above two inequality,
	\begin{align}\label{WJ3yr}
	&{\widetilde L^{{\alpha _1}, \ldots ,{\alpha _m}}}\notag\\
	\lesssim&\|b\|_{BMO}\nu_{\vec{\omega}}(B)^{m(1-\kappa)}\sum_{j=1}^\infty(j+1)\frac{1}{\nu_{\vec{\omega}}(2^{j+1}B)^m}\prod_{i=1}^m\omega_i\big(2^{j+1}B\big)\bigg\|\Phi\bigg(\frac{|f_i|}{\lambda}\bigg)\bigg\|_{L\log L(\omega_i),2^{j+1}B}\notag\\
	\lesssim&\|b\|_{BMO}\nu_{\vec{\omega}}(B)^{m(1-\kappa)}
	\sum_{j=1}^\infty(j+1)\frac{1}{\nu_{\vec{\omega}}(2^{j+1}B)^m}\prod_{i=1}^m\bigg\|\Phi\bigg(\frac{|f_i|}{\lambda}\bigg)\bigg\|_{(L\log L)^{1,\kappa}(\omega_i)}
	\prod_{i=1}^m \omega_i(2^{j+1}B)^{\kappa}\notag\\
	\lesssim&\|b\|_{BMO}\prod_{i=1}^m\bigg\|\Phi\bigg(\frac{|f_i|}{\lambda}\bigg)\bigg\|_{(L\log L)^{1,\kappa}(\omega_i)}\sum_{j=1}^\infty(j+1)\frac{\nu_{\vec{\omega}}(B)^{m(1-\kappa)}}{\nu_{\vec{\omega}}(2^{j+1}B)^{m(1-\kappa)}}\notag\\
	\lesssim&\|b\|_{BMO}\prod_{i=1}^m\bigg\|\Phi\bigg(\frac{|f_i|}{\lambda}\bigg)\bigg\|_{(L\log L)^{1,\kappa}(\omega_i)}\notag.
	\end{align}
	The proof of Theorem \ref{thm3} is finished.
\end{proof}

\end{document}